\documentclass[final]{siamltex}
\usepackage[includehead, includefoot,hmargin={3.4cm,3.4cm},vmargin={2.85cm, 2.8cm}]{geometry}
\usepackage[latin1]{inputenc}
\usepackage[T1]{fontenc}
\usepackage{amsmath}
\usepackage{amsfonts}
\usepackage{subfigure}
\usepackage{subfigmat}
\usepackage{algorithm,algorithmic}
\usepackage{amssymb} 
\usepackage{url}
\newcommand{\R}{\mathbb R}
\newcommand{\Z}{\mathbb Z}
\newcommand{\N}{\mathbb N}

  \usepackage{graphicx}
\newcommand{\dpar}[2]{\dfrac{\partial #1}{\partial #2}}

\newcommand{\E}{{\mathcal{ E}}}

\usepackage[usenames,dvipsnames]{color}
\usepackage{bm}
\def\beq{\begin{equation}}
\def\eeq{\end{equation}}
\def\esplit{\end{split}}
\def\beqalign{\begin{array}{rl}}
\def\eeqalign{\end{array}}
\DeclareMathOperator*{\argmin}{\arg\!\min}
\DeclareMathOperator*{\argmax}{\arg\!\max}
\makeatletter
\newcommand{\ALOOP}[1]{\ALC@it\algorithmicloop\ #1%
  \begin{ALC@loop}}
\newcommand{\ENDALOOP}{\end{ALC@loop}\ALC@it\algorithmicendloop}

\makeatother

\def\Abold{\mathbf{A}}

\def\Ebold{\mathbf{E}}

\def\Jbold{\mathbf{J}}

\def\Sbold{\mathbf{S}}

\def\Vbold{\mathbf{V}}
\def\Wbold{\mathbf{W}}

\def\Zbold{\mathbf{Z}}

\def\bbold{\mathbf{b}}

\def\fbold{\mathbf{f}}
\def\gbold{\mathbf{g}}

\def\mbold{\mathbf{m}}

\def\qbold{\mathbf{q}}
\def\rbold{\mathbf{r}}
\def\sbold{\mathbf{s}}
\def\tbold{\mathbf{t}}
\def\ubold{\mathbf{u}}
\def\vbold{\mathbf{v}}
\def\wbold{\mathbf{w}}

\def\ybold{\mathbf{y}}
\def\zbold{\mathbf{z}}
\def\alphabold{\boldsymbol{\alpha}}

\def\mubold{\boldsymbol{\mu}}
\def\taubold{\boldsymbol{\tau}}
\def\rhobold{\boldsymbol{\rho}}

\def\0bold{\boldsymbol{0}}
\def\0{\mathbf{\0}}

\newtheorem{remark}{remark}

\title{Robust Model Reduction Of Hyperbolic Problems by $L^1$-norm Minimization and Dictionary Approximation}
\author{R. Abgrall $^{(1)}$, D. Amsallem $^{(2)}$ and R. Crisovan $^{(1)}$\\
(1): Institut f\"ur Mathematik, Winterthurstrasse 190, \\CH 8057  Z\"urich, Switzerland \\
(2): Department of Aeronautics and Astronautics, 
496 Lomita Mall, \\Stanford University, Stanford, CA 94305-3035, USA}

\begin{document}
\maketitle


\begin{abstract}
We propose a novel model reduction approach for the approximation of non linear hyperbolic equations
in the scalar and the system cases. The approach relies on an offline computation of a dictionary of solutions together with an online $L^1$-norm minimization of the residual. It is  shown why this is a natural framework for hyperbolic problems and tested on nonlinear problems such as Burgers' equation and the one-dimensional Euler equations involving shocks and discontinuities. Efficient algorithms are presented for the computation of the $L^1$-norm minimizer, both in the cases of linear and nonlinear residuals. Results indicate that the method has the potential of being accurate when involving only very few modes,  generating physically acceptable, oscillation-free, solutions.
\end{abstract}

\begin{keywords} 
Reduced order models, non linear hyperbolic problems, $L^1$ minimisation.
\end{keywords}

\begin{AMS}
65M08, 65K99, 76L99 
\end{AMS}

\pagestyle{myheadings}
\thispagestyle{plain}
\markboth{R. Abgrall, D. Amsallem and R. Crisovan}{Robust Model Reduction Of Hyperbolic Problems by $L^1$-norm Minimization}


\section{Introduction}
Model reduction is becoming an essential tool to enable application requiring either real-time predictions or the evaluation of a large number of partial differential equations (PDE) based computational models. The first category encompasses optimal control~\cite{ito98,ly01} and model predictive control~\cite{hovland08,amsallem13:mpc}. Routine analysis and parameterized studies~\cite{amsallem10}, design optimization~\cite{legresley00,amsallem15:smo} and the quantification of uncertainty~\cite{buithanh08} are applications pertaining to the second category, to name just a few. In all of these applications, the large dimensionality associated with the discretized partial equations prevents their solution in real-time. Model reduction reduces that cost by restricting the solution to a subspace of the solution space. This subspace is usually described by a small number of reduced basis vectors. In turn, a projection step reduces the dimensionality of the system of discrete equations considered, enabling their fast solution. 

While the model reduction of elliptic and parabolic PDEs has been the subject of numerous studies~\cite{kunisch02,veroy03} and its theory is well understood~\cite{kunisch01,grepl05,rozza08}, reducing hyperbolic equations has proved to be much more challenging~\cite{barone09}. More specifically, moving waves and discontinuities such as shocks require a large number of basis vectors to accurately approximate these features~\cite{dahmen14}. This characterizes these problems as ones with large Kolmogorov $n$-widths~\cite{binev11}.

To circumvent this issue, approaches based on local bases~\cite{amsallem08,dihlmann11,amsallem12:localROB,maday12} reduce the Kolmogorov $n$-width by considering local subspaces. The locality can be characterized in parameters~\cite{amsallem08}, time~\cite{dihlmann11} or state-space~\cite{amsallem12:localROB}. In the present work, an approach based on dictionaries is considered~\cite{kaulmann13,brunton13}. More specifically, solutions corresponding to various time and parameter instances are collected and stored in such a dictionary. Each solution will then be considered as a reduced basis vector. In turn, localization in time and space can be easily enforced by only considering basis vectors corresponding to restricted subdomains of the time and parameters spaces. In addition to the reduction in number of basis vectors, this paper will demonstrate that a key advantage of a dictionary approach is a better approximation of states having sharp gradients and discontinuities. In particular, it will be demonstrated that avoiding basis truncation such as the one occurring in Proper Orthogonal Decomposition (POD)~\cite{sirovich87} or Non-Negative Matrix Factorization~\cite{balajewicz15} avoid Gibbs phenomenon.

In addition to the choice of reduced basis, a key ingredient in projection-based model reduction is the definition of the reduced system of equations. For symmetric systems such as those arising in elliptic and parabolic PDEs, Galerkin projection is the method of choice. For nonsymmetric systems, however, it has been shown that minimizing the $L^2$-norm of the residual is preferable for stability considerations~\cite{carlberg11,carlberg13}. In the present paper, model reduction based on the minimization of the $L^1$-norm of the residual is introduced and its advantage is demonstrated in conjunction with a dictionary approach for reducing problems with sharp gradient and shocks. More specifically, the present work demonstrates that combining a dictionary and $L^1$ minimization promotes sparsity in the choice of basis functions participating in the reduced-order solution and results in more accurate and physical reduced-order solutions.

This paper presents practical algorithms for performing $L^1$-norm minimization both in the linear and nonlinear cases. Furthermore, in order to achieve practical speed-up, another level of approximation, hyper-reduction~\cite{ryckelynck05,chaturantabut10,carlberg11,carlberg13} is required. As such, hyper-reduced versions of $L^1$-norm residual minimization are developed as well.

This paper is organized as follows: we first discuss the problem of interest: the approximation of the solution of non linear problems by reduced order models. In the second section, we  explain the role of $L^1$ minimization in this problem, then discuss some difficulties associated to this. Then we present in detail the algorithm we have developed, and provide an error estimate. The last section provides several numerical example that illustrate the behavior of our methods, on linear and non linear problems. A conclusion follows and we sketch some perspectives.

\section{Problem of interest}
In this work, high-dimensional models (HDM) arising from the space discretization of hyperbolic PDEs are considered. PDEs of the following type are considered
\begin{equation}\label{eq:PDE}
\displaystyle{\left\{ 
  \begin{array}{l l }
 \displaystyle{ \frac{\partial W}{\partial t}} + L(W,\mubold) &= f(t,\mubold),~x\in\Omega,~t\in[0,T], \\
B(W,\mubold) &= g(t,\mubold),~x\in\partial\Omega,~t\in[0,T],\\
W(x,t=0,\mubold) &= W_0(x,\mubold),~x\in\Omega
\end{array} \right.}
\end{equation}
$W\in\mathbb{R}^p$ is a scalar ($p=1$) or vector ($p>1$) field, $\Omega\subset\mathbb{R}^d$ is the domain of the equation $1\leq d\leq 3$ and $\partial\Omega$ the boundary of the domain. $L$ is a differential operator such as the Laplacian or the divergence of a flux and $B$ a boundary operator, $f$ and $g$ are volume and surface forces, respectively and  $\mubold\in\mathcal{P}\subset\mathbb{R}^{m}$ is a vector of $m$ parameters defining the system of interest.

The HDMs result from a finite differences approximation or finite volume formulation of the PDE~(\ref{eq:PDE}) under the following form
\begin{equation}\label{eq:ODE}
\displaystyle{\left\{ 
  \begin{array}{l l }
 \displaystyle{ \frac{d\wbold}{dt} + \fbold(\wbold(t),t,\mubold)} &= \gbold(t,\mubold),~t\in[0,T]\\
\wbold(t=0) &= \wbold_0,
\end{array} \right.}
\end{equation}
where $\wbold\in\mathbb{R}^N$ is the HDM state of large dimension $N$, $t$ denotes time. $\fbold(\cdot,\cdot)$ and $\gbold(\cdot)$ are nonlinear functions of their arguments. 

 In the remainder of this paper, the time and parameter variables are grouped together, unless explicitly stated, as a variable $\taubold = [t;\mubold]$. Hence, the HDM state is parameterized as
\begin{equation}
\wbold(\taubold) = \wbold(t,\mubold).
\end{equation}
In practice, the ODE~(\ref{eq:ODE}) is discretized in time using a time discretization $t_0=0<t_1<\cdots<t_{N_t}=T$. Explicit and implicit time-discretization techniques are used in the present paper, resulting in a sequence of nonlinear systems of equations of large dimension $N$
\begin{equation}\label{eq:HDMres}
\rbold^n(\wbold) = \boldsymbol{0},~n=1,\cdots,N_t,
\end{equation}
where $\rbold^n = [r^n_1,\cdots,r^n_N]^T$. We give several examples later in the text, note that the residual $\rbold^n$ will depend of several time instances of the solution for unsteady problems, for example $\wbold^n$ and $\wbold^{n-1}$ in the simplest case. Steady problems can also be written in the form $\rbold(\wbold)=\boldsymbol{0}$. Since our goal is to draw a proof of concept, we will only take two generic examples later in the text. 

The goal of model reduction is to approximate the high-dimensional system~(\ref{eq:HDMres}) using a much smaller number of variables while retaining accuracy of the solution. For that purpose, projection-based model reduction techniques approximate the state $\wbold(\tau)$ in a subspace of $\mathbb{R}^N$ using a reduced-order basis (ROB) $\Vbold = [\vbold_1,\cdots,\vbold_k]\in\mathbb{R}^{N\times k}$. The state is then approximated as
\begin{equation}\label{eq:subsApprox}
\wbold(\tau) \approx \Vbold\qbold(\tau) = \sum_{i=1}^k \vbold_i q_i(\tau)
\end{equation}
where $\qbold(\tau) = [q_1(\tau),\cdots,q_k(\tau)]^T$ denotes the vector of $k$ reduced coordinates. Substituting the subspace approximation~(\ref{eq:subsApprox}) into~(\ref{eq:HDMres}) usually results in a non-zero residual of dimension $N$
\begin{equation}\label{eq:nnzRes}
\rbold^n(\Vbold\qbold) \approx \boldsymbol{0}.
\end{equation}
Two common approaches result in the definition of a reduced system of equations:
\begin{itemize}
\item Galerkin projection enforces the orthogonality of the residual to the ROB $\Vbold$ as
\begin{equation}
\Vbold^T\rbold^n(\Vbold\qbold) = \boldsymbol{0},~n=1,\cdots,N_t.
\end{equation}
This defines a set of $k$ nonlinear equations in terms of $k$ unknowns which can be solved by Newton-Raphson's method.
\item Residual minimization approaches~\cite{legresley00,buithanh08,carlberg11,amsallem12:localROB,carlberg13} minimize the residual in the $L^2$-norm sense
\begin{equation}
\min_{\qbold} \| \rbold^n(\Vbold\qbold)\|_2^2=\sum_{i=1}^N \left(r_i^n(\Vbold\qbold)\right)^2,~n=1,\cdots,N_t.
\end{equation} 
In practice, this nonlinear least-squares problem can be solved using Gauss-Newton or Levenberg-Marquardt iterations~\cite{nocedalbook}. In Section~\ref{sec:L1}, alternative residual minimization approaches based on $L^1$-norm minimization which are more appropriate for the reduction of hyperbolic problems will be proposed.
\end{itemize}
\section{Dictionary approach}\label{sec:dico}
Projection-based model reduction techniques~\cite{sirovich87,rozza08,balajewicz15} based on snapshots pre-compute solutions of the HDM for specific values of the vector $\taubold=[t;\mubold]$. These snapshots are gathered in a snapshot matrix
\begin{equation}
\Sbold = [\wbold(\taubold_1),\cdots,\wbold(\taubold_{N_s})].
\end{equation}
Two approaches for compressing the snapshot matrix are described as follows:
\begin{itemize}
\item Proper Orthogonal Decomposition~\cite{sirovich87} computes an optimal reduced-order basis of dimension $k$ that minimizes the projection error of the snapshots onto the basis.
\item Balanced POD~\cite{willcox01}, applicable to linear systems only, also takes into account snapshots of the dual system to construct the reduced basis for the primal and dual systems.
\item Non-negative matrix factorization~\cite{lee99} was recently applied to construct a non-negative reduced-order basis based on snapshots with positive entries in the context of contact problems~\cite{balajewicz15}. The reduced basis minimizes the positive reconstruction of the snapshots.
\end{itemize}
All three approaches perform a compression of the information contained in the snapshot matrix $\Sbold$. More specifically, the $N_s$ vectors contained in $\Sbold$ are compressed, leading to a reduced-order basis of dimension $k\leq N_s$.

In the present paper, an approach based on a dictionary of solutions is preferred as it does not incur any loss of information by compression. As such, the vectors $\{\vbold_i\}_{i=1}^k$ in the reduced basis are solutions of the HDM:
\begin{equation}
\vbold_i = \wbold(\taubold_i),~i=1\cdots,k.
\end{equation}
The solution of the HDM will then be approximated as
\begin{equation}
\wbold(\taubold) \approx \sum_{i=1}^k \wbold(\taubold_i)\qbold(\taubold).
\end{equation}

In the present case, since the HDM is of very large dimension, over-complete dictionaries, as used in compressed sensing~\cite{candes06,donoho06} and for which $k\geq N$ will not be considered.

\section{$L^1$-norm residual minimization}\label{sec:L1}
In the present paper, model reduction based on $L^1$-norm residual minimization is introduced to reduce the dimensionality of hyperbolic equations as an alternative to Galerkin projection and $L^2$-norm minimization. Motivations for the use of the $L^1$-norm are provided in Section~\ref{ssec:L1motiv}. Model reduction based on $L^1$-norm minimization is introduced in Section~\ref{ssec:L1MOR} together with practical numerical procedure for their computation in Section~\ref{ssec:L1Alg}. 
 
\subsection{Motivations}\label{ssec:L1motiv}
Minimizing the $L^1$-norm of the residual is known to lead to regressions that are much more robust to outliers~\cite{boyd04}. In the context of hyperbolic systems, the work of Guermond et al. on Hamilton Jacobi equations and transport problems \cite{guermond08,guermond09} has shown, at least experimentally, that the numerical solution can retain an excellent non-oscillatory behavior by minimizing the $L^1$-norm of the PDE residual.  For completeness, the motivation for $L^1$-norm minimization is justified as follows for the problem
\begin{equation}
\label{eq:2}
\dpar{W}{t}+\text{ div } {F}(W)=0
\end{equation}
defined on $\Omega \subset \R^d$ and for $t> 0$.  The solution $W$ belongs here to $\R^p$, so that $F=(F_1, \ldots, F_p)^T$.
The weak form of the equation is: for any $\varphi \in \left[C^1(\Omega)\right]^p$ and with compact support:
\begin{equation}\int_\Omega \varphi(x,t)\bigg ( \dpar{W}{t}+\text{ div } F(W) \bigg )dx=0.
\end{equation}
Integrating by parts yields
\begin{equation}
\int_\Omega \dpar{\varphi}{t} W dx+\int_\Omega \nabla \varphi \cdot F(W) dx=0.
\end{equation}
Restricting to the set of test functions $\left\{\varphi\in \left[C^1(\Omega)\right]^p, ||\varphi||_{\infty}\leq 1\right\}$,  $W$ is a solution if:
\begin{equation}\label{remi:eq:1}
\sup\limits_{\{\varphi\in \left[C^1(\Omega)\right]^p, ||\varphi||_{\infty}\leq 1\}}\Bigg (\int_\Omega \dpar{\varphi}{t} W dx+\int_\Omega \nabla \varphi \cdot F(W) dx\Bigg )=0.
\end{equation}
Remember that for  any function $g\in L^1(\R^d)$, the total variation is defined as
\begin{equation}
TV(g)=\sup\limits_{\varphi\in C^1_0(\R^d)\cap L^\infty(\R^d), ||\varphi||_{\infty}\leq 1} \left\{\int_{\R^d} \nabla \varphi (x)\cdot g(x) dx\right\},
\end{equation}
and if, in addition, $g\in C^1(\R^d)$, then $TV(g)=\int_{\R^d}||\nabla g||dx=||\nabla g||_{L^1(\R^d)}.$
This shows that, defining the space-time flux $\mathcal{F}=(W,F)$, $W$ is a weak solution if and only if the (space-time)total variation of $\mathcal{F}$ vanishes, that is
\begin{equation}
\label{rem:justificationTV}
 TV\big ( \mathcal{F}(W)) =0.
 \end{equation}
 In other works, one can look for $W$ as a function of $L^1\cap L^\infty$ such that $W$ minimizes
 $TV\big ( \mathcal{F}(V))$ over $V\in L^1\cap L^\infty$, i.e.
 \begin{equation}
 \label{remi:argmin}W=\text{argmin}\{  TV\big ( \mathcal{F}(V)), V\in L^1\cap L^\infty\}.\end{equation}
 This does not garanty uniqueness (and thus there is some abuse of language in this setting), since the entropy conditions are not encoded into this formulation.
 However, \eqref{remi:argmin} indicates that a natural setting is to minimize the $L^1$ norm of the space-time divergence of the space-time flux $\mathcal{F}$.
 
 \bigskip

How does it translates in the discrete setting? For simplicity, we only mention the case of explicit schemes. The case  of implicit ones is done similarly, we discuss later the solution procedure.

The following useful classical result is mentioned.
Consider $\{x_i\}_{i\in \Z}$ a strictly increasing sequence in $\R$ and $x_{i+1/2}=\frac{x_i+x_{i+1}}{2}$. Assuming that $\R=\cup_{i\in \Z}[x_{i-1/2},x_{i+1/2}[$ and considering $g$ defined by, for any $i\in \Z$, 
\begin{equation}
g(x)=g_i \text{ if } x\in [x_{i-1/2},x_{i+1/2}[,
\end{equation}
then
\begin{equation}
TV(g)=\sum_{i\in \Z} |g_{i+1}-g_i|.
\end{equation}

Now, instead of having the exact solution, consider an approximation procedure that enables, from $\wbold^{n}\approx W(~.~, t_n)$, to compute $\wbold^{n+1}\approx W(~.~, t_{n+1})$,
say $\mathcal{L}(\wbold^n, \wbold^{n+1})$.

For instance, assume that  we have a finite volume method and $d=1$: for any grid point $i\in\{1,\cdots,N\}$,
\begin{equation}
\left[\mathcal{L}(\wbold^n, \wbold^{n+1})\right]_i= \Delta x(\wbold_i^{n+1}-\wbold_i^n)+ {\Delta t}\big (\fbold_{i+1/2}(\wbold^n)-\fbold_{i-1/2}(\wbold^n)\big ).
\end{equation}
Here $\fbold_{j+1/2}$ is any numerical flux at the cell interface $x_{j+1/2}$, see for example \cite{toro} for the classical examples. 

One way to evaluate $\wbold^{n+1}$ is to minimize the total variation, i.e.
\begin{equation}TV(\mathcal{L})=\sum_{i\in \mathcal{I}} \Big |\Delta x (\wbold_i^{n+1}-\wbold_i^n)+ {\Delta t}\big(\fbold_{i+1/2}(\wbold^n)-\fbold_{i-1/2}(\wbold^n)\big )\Big |,
\end{equation}
leading to
\begin{equation}
\wbold^{n+1}= \argmin _{\vbold \text{ piecewise constant }} \sum_{i\in \mathcal{I}} \bigg |\Delta x (\vbold_i-\wbold_i^n)+ \Delta t(\fbold_{i+1/2}(\wbold^n)-\fbold_{i-1/2}(\wbold^n))\bigg |.
\end{equation}
Clearly, if $\mathcal{I}$ is equal to the set of grid points, the solution is given by 
\begin{equation}
\wbold_i^{n+1}=\wbold_i^n -\dfrac{\Delta t}{\Delta x} \bigg (\fbold_{i+1/2}(\wbold^n)-\fbold_{i-1/2}(\wbold^n)\bigg ).
\end{equation}

\subsection{Model reduction by $L^1$-norm  minimization}\label{ssec:L1MOR}
As an alternative to Galerkin projection and residual minimization in the least-squares sense, a reduced system of equation is here obtained by minimizing, at each time step $n=1,\cdots,N_t$, the $L^1$-norm of the residual vector as
\begin{equation}\label{eq:L1minres}
\min_{\qbold} \| \rbold^n(\Vbold\qbold)\|_1 = \sum_{i=1}^n \left| r^n_i(\Vbold\qbold)\right|,~n=1,\cdots,N_t
\end{equation}

or as 

\begin{equation}\label{eq:L1minres_convex_hull}
\min_{\qbold} \| \rbold^n(\Vbold\qbold)\|_1 \quad subject \quad to\quad \textbf{1}^{T}q=1,q\geqslant 0, ~n=1,\cdots,N_t,
\end{equation}
where we project on the convex envelop of the dictionary (see Section \ref{sec:Err_Uniq} for more details). There are at least two difficulties associated with minimizing the $L^1$-norm. A first one  is is that the $L^1$ norm is not strictly convex, so that the uniqueness is not guaranteed. This difficulty is taken into account in the solution procedure by adding a strictly convex penalization term, for example a $L^2$ constraint. The second difficulty associated to $L^1$ is  its non-differentiability at zero. To circumvent this issue, the Huber function~\cite{huber11}, defined as follows can be introduced:
\begin{equation}\label{eq:Huberfun}
\phi_M(x) = \left\{  
 \begin{array}{l l }
 x^2 & \text{if}~ |x|\leq M \\
M(2|x|-M) & \text{otherwise} ,
\end{array}
 \right.
\end{equation}
Then, the sequence of reduced systems of equations based on the Huber function is
\begin{equation}\label{eq:Huberminres}
\min_{\qbold}  \sum_{i=1}^n \phi_M\left( r^n_i(\Vbold\qbold)\right),~n=1,\cdots,N_t.
\end{equation}

The Huber function $\phi_M$ behaves as a parabola close to $x=0$ and as the $L^1$-norm for large values of $x$. It is continuously differentiable on $\mathbb{R}$ ($\phi_M\in\mathcal{C}^1(\mathbb{R})$). It is also used in regressions as a loss function due to its non-sensitivity to outliers. In the present work, it will be used as a continuously differentiable alternative to the $L^1$-norm.

\begin{figure}[h!]
\begin{center}
{\includegraphics[width=0.75\textwidth,clip=]{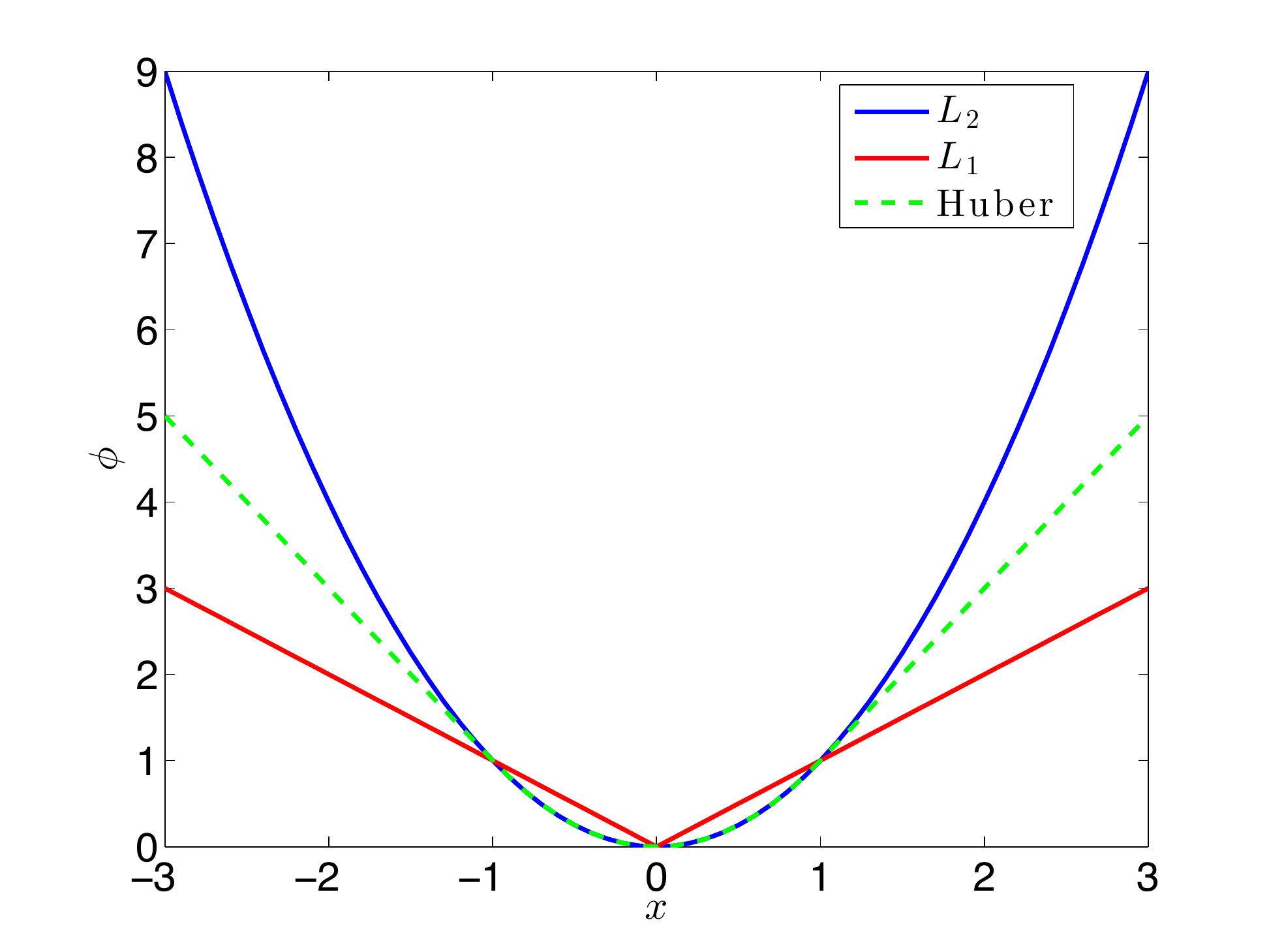}}
\end{center}
\caption{\label{fig:NormsComp}Comparison of the $L^2$, $L^1$ and Huber function ($M=1$) norms}
\end{figure}

Figure~\ref{fig:NormsComp} compares, in the scalar case, the $L^2$ and $L^1$-norms to the norm based on the Huber function for the particular case $M=1$.
Practical algorithm for solving the systems of equations~(\ref{eq:L1minres}) and~(\ref{eq:Huberminres}), both in the case of linear and nonlinear residual functions are presented in the following section.

\subsection{Algorithms}\label{ssec:L1Alg}
A classical solution to minimizing a linear residual vector in the $L^1$-norm is by recasting the problem as a linear program (LP). More specifically, assuming that the residual is linear $\rbold^n(\Vbold\qbold) = \Abold^n\Vbold\qbold + \bbold^n$ with $\Abold^n\in\mathbb{R}^{N\times N}$ and $\bbold^n\in\mathbb{R}^N$, a solution to~(\ref{eq:L1minres}) is given by the solution $\qbold$ of the LP
\begin{equation}\label{eq:LP}
\begin{aligned}
 \min_{\qbold,\sbold,\tbold}~ \boldsymbol{1}^T(\sbold+\tbold)&\\
\text{s.t.}~\Abold^n\Vbold\qbold + \bbold^n -\sbold + \tbold  &=   \boldsymbol{0}\\
            \sbold &\geq \boldsymbol{0}\\
            \tbold &\geq \boldsymbol{0}.
\end{aligned}
\end{equation}
Unfortunately, this LP involves $k+2N$ variables and $3N$ constraints, including $N$ equality constraints, rendering this approach intractable in the case of model reduction. 

Alternatively, the $L^1$-norm minimization problem can be solved by Iteratively Reweighted Least Squares (IRLS)~\cite{daubechies08}.
This approach proceeds iteratively by solving a sequence of weighted least-squares problem. An advantage of this approach is that its implementation can rely entirely on existing least-squares solvers. Furthermore, its complexity is similar to that of the $L^2$-norm minimization problem. The procedure is presented in Algorithm~\ref{alg:L1IRLS} in the case of a nonlinear residual vector. At each iteration $l$, a weighted least-squares problem is solved, where the weight depend on the current value of the residual vector $\rbold^l$ as follows:
$\Wbold^l = {\bf{diag}}\left({|r^l_i|}^{-\frac{1}{2}}\right)$

\begin{algorithm}
\fontsize{10pt}{10pt}\selectfont
\caption{$L^1$-norm minimization  by  Iteratively Reweighted Least-Squares} 
\begin{algorithmic}[1]
{\REQUIRE Residual function $\rbold(\cdot)$ and associated Jacobian $\Jbold(\cdot)$, reduced basis $\Vbold$, initial guess $\qbold^0$, tolerance for convergence $\epsilon$
\ENSURE Solution $\qbold$
\STATE $l=0$
\STATE Compute $\rbold^0 =\rbold(\Vbold\qbold^0)$ and $\Zbold^0 =\Jbold(\Vbold\qbold^0)\Vbold$
 \WHILE { $l=0$ OR $ \| \Delta \qbold^{l-1}|_1 >  \epsilon (1+ \| \qbold^{l-1}\|_1)$} 
 \STATE Compute the weights  $\Wbold^l = \text{diag}\left({|r^l_i|}^{-\frac{1}{2}}\right)$
 \STATE Solve the weighted least-squares problem
  \[ \Delta\qbold^l =  \argmin_{\ybold}\| \Wbold^l\Zbold^l \ybold + \Wbold^l\rbold^l\|_2^2 \]
 \STATE $\qbold^{l+1} = \qbold^l + \Delta \qbold^l$
 \STATE Compute $\rbold^{l+1} =\rbold(\Vbold\qbold^{l+1})$ and $\Zbold^{l+1} =\Jbold(\Vbold\qbold^{l+1})\Vbold$
   \STATE $l=l+1$
 \ENDWHILE
 \STATE $\qbold = \qbold^l$
 }
 \end{algorithmic}\label{alg:L1IRLS}
\end{algorithm}

Similarly, minimization of the Huber function can also be done by an IRLS procedure, as described in Algorithm~\ref{alg:HuberIRLS}. The procedure only differs from its $L^1$-norm counterpart by the choice of weights. In the present work, the following choice of weights is proposed for a given residual vector $\rbold^l$
\begin{equation}
\Wbold^l = {\bf{diag}}\left(\delta(|r_i^l|< M) + M{|r^l_i|}^{-\frac{1}{2}}\delta(|r_i^l|\geq M) \right).
\end{equation}
Furthermore, the parameter $M$ is computed as
\begin{equation}\label{eq:M}
M=\epsilon_2\max(1,\max(|r^l_i|))
\end{equation}
with $\epsilon_2=10^{-6}$, as it has been found to be a robust choice across different applications.
\begin{algorithm}
\fontsize{10pt}{10pt}\selectfont
\caption{Huber function minimization by Iteratively Reweighted Least-Squares} 
\begin{algorithmic}[1]
{\REQUIRE Residual function $\rbold(\cdot)$ and associated Jacobian $\Jbold(\cdot)$, reduced basis $\Vbold$, initial guess $\qbold^0$, tolerance for convergence $\epsilon$
\ENSURE Solution $\zbold$
\STATE $l=0$
\STATE Compute $\rbold^0 =\rbold(\Vbold\qbold^0)$ and $\Zbold^0 =\Jbold(\Vbold\qbold^0)\Vbold$
 \WHILE { $l=0$ OR $ \| \Delta \qbold^{l-1}|_1 >  \epsilon (1+ \| \qbold^{l-1}\|_1)$} 
 \STATE Compute the weights $\Wbold^l = {\bf{diag}}\left(\delta(|r_i^l|< M) + M{|r^l_i|}^{-\frac{1}{2}}\delta(|r_i^l|\geq M) \right)$
 \STATE Let $M=\epsilon_2\max(1,\max(|r^l_i|))$
 \STATE Solve the weighted least-squares problem
  \[ \Delta\qbold^l =  \argmin_{\ybold}\| \Wbold^l\Zbold^l \ybold + \Wbold^l\rbold^l\|_2^2 \]
 \STATE $\qbold^{l+1} = \qbold^l + \Delta \qbold^l$
 \STATE Compute $\rbold^{l+1} =\rbold(\Vbold\qbold^{l+1})$ and $\Zbold^{l+1} =\Jbold(\Vbold\qbold^{l+1})\Vbold$
   \STATE $l=l+1$
 \ENDWHILE
 \STATE $\qbold = \qbold^l$
 }
 \end{algorithmic}\label{alg:HuberIRLS}
\end{algorithm}

\section{Model Reduction by $L^1$-norm minimization and dictionaries}
\subsection{Procedure}
In this section, model reduction based on minimizing the residual in the $L^1$-norm is combined with the dictionary approach presented in Section~\ref{sec:dico}. 

A potential issue with using a dictionary, as opposed as a reduced basis, is the fact that the dictionary may be rank-deficient. One option to address this issue is to perform a Gramm-Schmidt orthogonalization or  a rank-revealing QR factorization. A drawback of that approach is that dictionary members are then linearly combined.  Alternatively, a regularization term is here added to the minimization functionals to ensure a system with full rank and a unique solution as follows (see Section~\ref{sec:Err_Uniq} for more details).
\begin{itemize}
\item For $L^1$-norm minimization, the functional becomes
\begin{equation}\label{eq:L1minresReg}
\min_{\qbold} \| \rbold^n(\Vbold\qbold)\|_1 + \eta \| \qbold\|_2^2 = \min_{\qbold} \sum_{i=1}^N \left| r^n_i(\Vbold\qbold)\right| + \eta\sum_{j=1}^k q_j^2,~n=1,\cdots,N_t.
\end{equation}
\item For Huber function minimization, the functional becomes
\begin{equation}\label{eq:HuberminresReg2}
\min_{\qbold}  \sum_{i=1}^N \phi_M\left( r^n_i(\Vbold\qbold)\right) + \eta \| \qbold\|_2^2,~n=1,\cdots,N_t.
\end{equation}

\end{itemize}
\subsection{Selection of the dictionary members}
With large dictionaries, it may be computationally expensive to consider all dictionary members as potential basis vectors for a given value $\taubold^\star=(t^\star,\mubold^\star)$ of time and parameters. Instead, a local dictionary approach can be consider by restricting the dictionary members considered for $\taubold^\star$ in a neighborhood of the time and/or parameter domains. The local dictionary is then defined as
\begin{equation}
\Vbold(\taubold^\star) = [\wbold(\taubold_1),\cdots,\wbold(\taubold_r)]~,\taubold_i\in \mathcal{T}(t^\star)\times \mathcal{P}(\mubold^\star)\subset [0,T] \times \mathcal{P}.
\end{equation}
In such a dictionary approach, restricting dictionary members in the time and/or parameter domains is straightforward, unlike the case of pre-computed reduced bases, for which an a priori partitioning of the parameter domains are necessary~\cite{dihlmann11,amsallem12:localROB}.

\subsection{Training by greedy sampling}\label{sec:training}
An essential step in the construction of a parametric ROM is the selection of the sampled snapshots in the time and parameter domains. Greedy approaches~\cite{grepl05,buithanh08,pdt15} proceed by iteratively selected the location in the parameter space where the error between the HDM and the ROM is the largest. As computing the error requires the expensive solution of the HDM, cheaper error indicators are used instead. In the present work, the cumulated $L^1$-norm of the residual vector corresponding to the ROM solution is used as error indicator:
\begin{equation}\label{eq:errorIndicL1}
\mathcal{I}(\mubold) = \sum_{n=1}^{N_t} \| \rbold^n(\Vbold\qbold^n(\mubold))\|_1,
\end{equation}
where $\qbold^n(\mubold)$ denotes the ROM solution for the parameter $\mubold$ at time iteration $t^n$. The greedy procedure is detailed in Algorithm~\ref{alg:greedyL1}.

\begin{algorithm}
\fontsize{10pt}{10pt}\selectfont
\caption{Greedy sampling of the parameter space}
\begin{algorithmic}[1]
{\REQUIRE Residual function $\rbold(\cdot)$, tolerance for convergence $\epsilon$, candidate parameter set $\mathcal{C} = \{ \mubold^{(i)}\}_{i=1}^{N_c}$
\ENSURE Dictionary $\Vbold$
\STATE Randomly chose an initial sample parameter $\mubold_0\in\mathcal{C}$ and compute the associated HDM solution $\{\wbold^n(\mubold_0)\}_{n=1}^{N_t}$
\STATE Construct an initial dictionary $\Vbold = \{\wbold^n(\mubold_0)\}_{n=1}^{N_t}$
\FOR {$i_c=1,\cdots,N_c$}
\STATE Solve for the ROM solution $\{\qbold^n(\mubold^{(i_c)}\}_{n=1}^{N_t}$ and evaluate the error indicator $\mathcal{I}(\mubold^{(i_c)})$
\ENDFOR
\STATE $j=1$
 \WHILE { $\max_{i_c=1,\cdots,N_c} \mathcal{I}(\mubold^{(i_c)}) > \epsilon$} 
 \STATE Select $\mubold_j = \argmax_{i_c=1,\cdots,N_c} \mathcal{I}(\mubold^{(i_c)})$
 \STATE Compute the associated HDM solution $\{\wbold^n(\mubold_j)\}_{n=1}^{N_t}$
 \STATE Update the dictionary $\Vbold = \Vbold \bigcup \{\wbold^n(\mubold_j)\}_{n=1}^{N_t}$
\FOR {$i_c=1,\cdots,N_c$}
\STATE Solve for the ROM solution $\{\qbold^n(\mubold^{(i_c)}\}_{n=1}^{N_t}$ and evaluate the error indicator $\mathcal{I}(\mubold^{(i_c)})$
\ENDFOR
   \STATE $j=j+1$
 \ENDWHILE
 }
 \end{algorithmic}\label{alg:greedyL1}
\end{algorithm}
\section{Error estimation}\label{sec:Err_Uniq}
In this section, we provide an error estimate (in the scalar case) between the solution obtained by projection on the dictionary and the initial error. This error estimate is  another way to justify the method, in particular the greedy algorithm used for the dictionary selection. These estimates are provided in a simple setting: we consider monotone scheme. In this section, we first precise the setting, then give a natural condition on the dictionary for obtaining this error estimate and then state and prove it.

\subsection{Scheme setting} Consider the scalar conservation law equations with the initial condition:
\begin{equation} \label{problem}
\begin{split}
\dpar{u}{t}+\dpar{f(u)}{x} &= 0, \quad x\in \R, t>0\\
u(x,0) &= u_0(x), \quad x\in \R. 
\end{split}
\end{equation}
After discretizing, we assume that the scheme writes,  for $u:=(u_j)_{j\in \Z}$, 
\begin{equation}\label{remi:schema}
u^{n+1}=S(u^n, \lambda)
\end{equation} 
with  $\lambda=\Delta t/\Delta x$ and  the initial condition
\begin{equation}
\label{eq:1:2}
u_j^0=\text{given}.
\end{equation}
 We assume that the operator $S$ is  monotone for $\lambda \in [0,b[$, $b>0$, i.e. if for any sequence $u$ and $v$ bounded for the $L^1$ or $L^\infty$ norms  with $j\in \Z$, $u_j\leq v_j$, then $S(u,\lambda)_j\leq S(v,\lambda)_j$. Let $L^1$ and $L^\infty$ norms are generically denoted by $||~.~||$.

An example is given by the scheme
\label{eq:1}
\begin{equation}
\label{eq:1:1}
S(u)_j=u_j-\lambda \big ( \hat{f}(u_{j+1},u_j)-\hat{f}(u_j,u_{j-1})\big )
\end{equation}
where we assume that the numerical flux $\hat{f}(a,b)$ is monotone, i.e. increasing with respect to the first variable and decreasing w.r.t the second one. $S$ is monotone under a CFL like condition. Another example is given by the implicit scheme, where $v=S(u)$ is defined as the solution of
\begin{equation}
\label{eq:implicit}
v_j=u_j-\lambda \big ( \hat{f}(v_{j+1},v_j)-\hat{f}(v_j,v_{j-1})\big )
\end{equation} which is unconditionaly  monotone.

Thanks to Crandall-Tartar lemma (see \cite{God} for example), we know that for any $u$ and $v$, then 
$$||S(u,\lambda)-S(v,\lambda)||\leq ||u-v||$$
 in the  $L^1$ norm. The same is true in the $L^\infty$ norm.

\subsection{A standard result of convex programming}
If $b\in \R^n$ and $A\in M_{n,p}(\R)$ is a matrix with $n$ lines and $p$ columns, it is known that
the problem: find $x_{min}\in \R^p$ such that
$$||Ax_{min}-b||_1=\min_{y\in \R^p}||Ay-b||_1$$
has a unique solution, denoted by $p(b)$. In order to see this, we write   the minimisation problem as
$$\sum\limits_{i=1}^n |a_i^T x-b_i|,$$
where $a_i$ is the $i$-th line of $A$ and $b_i$ the $i$-th component of $b$. Hence it can be rewritten as the
minimisation of
$$\sum\limits_{i=1}^n u_i$$
under the constraints
$$
u_i\geq 0, \text{ and }
-u_i\leq a_i^T x-b_i\leq u_i
$$
i.e. with standard notations
$$
Ax-u\leq b, \text{ and }
-Ax-u\leq -b.
$$
Thus, this can be rewritten as the minimisation of $c^T X$ where $X= ( x_1,\ldots x_p, u_1, \ldots u_n)^T\in \R^{p+n}$, $c=(c_i)_{i=1, p+n}$ where $c_i=0$ for $i=1, \ldots p$ and $c_i=1$ for $i\geq p+1$ under the constraints
$$G X-B\leq 0$$
with 
$$G=\begin{pmatrix}
A & -\text{Id}_{n\times n}\\
-A &\text{Id}_{n\times n}
\end{pmatrix}\in M_{2n, p+n}(\R), \qquad B=(b, -b)^T\in \R^{2n}
$$

This is a standard problem of convex programming that has a unique solution provided that the 
rank of $G$ is larger than $p+n$, i.e. in this case if the rank of $A$ is larger than $p$.

\subsection{Error estimate}

We collect and store in a dictionary the solutions $\{u^{n}(\mu_{i})\}_{i}$ of the problem \eqref{problem} which correspond to various time and parameter instances and where the initial conditions are defined for the parameters $\{\mu_i\}_{i=1, \ldots , m}\in \mathcal{P}$.  The matrix $A$ is the matrix which $i$-th column corresponds to the discrete values that represent the $i$-th element of the dictionary at time $t_n$. Our main assumption is that the rank of this matrix is  equal to the number of element of the dictionary, for any time. This enables to define a projection  operator $p^n$ for any time $t_n$, by solving the minimization problem: knowing that $u^n_{\bm{\mu}} \in \underset{\mu_{i} \in \mathcal{P}}{span}(\{u^{n}(\mu_i)\}_{i})$, find $u^{n+1}_{\bm{\mu}} \in \underset{\mu_{i} \in \mathcal{P}}{span}(\{u^{n+1}(\mu_i)\}_{i})$ such that:

$$u^{n+1}_{\bm{\mu}}=\underset{\bm{\mu} \in \mathcal{P}}{argmin}\big \{ ||v_{\bm{\mu}}-S(u^{n}(\bm{\mu}), \lambda)||_1 : v_{\bm{\mu}} \in \underset{\bm{\mu} \in \mathcal{P}}{span}(\{u^{n+1}(\bm{\mu})\})\big \}=p^n(S(u^{n}(\bm{\mu}),\lambda).$$ 

We have  immediately the following  estimate:
\begin{equation}
\begin{split}
||p^n(S(u^{n}(\bm{\mu}),\lambda))-S(u^{n}(\bm{\mu}),\lambda)||_{1} &= \min_{ v_{\bm{\mu}} \in \underset{\bm{\mu} \in \mathcal{P}}{span}(\{u^{n+1}(\bm{\mu})\})}
 ||v_{\bm{\mu}}-S(u^{n}(\bm{\mu}),\lambda)||_{1}\\
& \leq \min_{\bm{\mu} \in \mathcal{P}} ||u_{\bm{\mu}}^{n+1}-S(u^{n}(\bm{\mu}),\lambda)||_{1} \\
&= \min_{\bm{\mu} \in \mathcal{P}} ||S(u_{\bm{\mu}}^{n},\lambda))-S(u^{n}(\bm{\mu}),\lambda)||_{1} \\
&
\leq \min_{\bm{\mu} \in \mathcal{P}} ||u_{\bm{\mu}}^{n}-u^{n}(\bm{\mu})||_{1} 
\end{split}
\end{equation} 
provided $\lambda$ enables to fulfill the monotonicity property for all the elements of the dictionary.
In the end, we get by induction the following estimate:
$$||p^n(S(u^{n}(\bm{\mu}),\lambda))-S(u^{n}(\bm{\mu}),\lambda)||_1\leq \min_{\bm{\mu} \in \mathcal{P}} ||u^{0}_{\bm{\mu}}-u^{0}(\bm{\mu})||_1
$$
where $u^0$ is the initial condition.\\

Next, we consider the case associated with the minimization problem (\ref{eq:L1minres_convex_hull}). If $\lambda_i\geq 0$ and $\sum_{i=1}^{|\mathcal{P}|} \lambda_i=1$, we obtain a sharp error estimate of type:\\

\begin{equation*}
\begin{split}
||p^n(S(u^{n}(\bm{\mu}),\lambda))-u^{n+1}(\bm{\mu})||_1&=||\sum_{i=1}^{|\mathcal{P}|} \lambda_i u_{\bm{\mu}_i}^{n+1}-u^{n+1}(\bm{\mu})||_1\\
&\leq \sum_{i=1}^{|\mathcal{P}|} |\lambda_i|\; ||u_{\bm{\mu}_i}^{n+1}-u^{n+1}(\bm{\mu})||_1\\
&\leq \big (\sum_{i=1}^{|\mathcal{P}|} |\lambda_i|\big ) \max_{\bm{\mu}\in \mathcal{P}}||u_{\bm{\mu}}^{n+1}-u^{n+1}(\bm{\mu})||_1\\
& =\max_{\bm{\mu}\in \mathcal{P}}||u_{\bm{\mu}}^{n+1}-u^{n+1}(\bm{\mu})||_1\\
& \leq\max_{\bm{\mu}\in \mathcal{P}} ||u_{\bm{\mu}}^{0}-u^{0}(\bm{\mu})||_1.
\end{split}
\end{equation*}


We have shown the following result:\\

\begin{proposition}
If $S(~.~,\lambda)$ is monotone for $\lambda\in [0, b[$ and if the dictionary $\{u^{n}(\mu_{i})\}_{i, \mu_i\in \mathcal{P}, n\geq 0}$ is of rank $|\mathcal{P}|$ for any $n\in \N$, then the reduced solution $u^{n+1}_{\bm{\mu}}=p^n(S(u^{n}(\bm{\mu}),\lambda))$  at time $t_{n+1}$ satisfies:
$$||p^n(S(u^{n}(\bm{\mu}),\lambda))-S(u^{n}(\bm{\mu}),\lambda)||_1\leq \min_{\bm{\mu} \in \mathcal{P}} ||u_{\bm{\mu}}^{0}-u^{0}(\bm{\mu})||_1.
$$
\end{proposition}

\begin{remark}
In practice, the rank condition is not always met. This is the case of the  problems of sections \ref{sec:burgers} and \ref{sec:euler}, at least for the initial time, since we have constant states initial. In order to meet this ran condition, we randomly perturb the component of the dictionary by a quantity $\varepsilon$ close to machine zero. Since the operator is monotone, this amounts to modify the error estimate by adding a constant that bound the perturbation,
$$||p^n(S(u^{n}(\bm{\mu}),\lambda))-S(u^{n}(\bm{\mu}),\lambda)||_1\leq \min_{\bm{\mu} \in \mathcal{P}} ||u_{\bm{\mu}}^{0}-u^{0}(\bm{\mu})||_1+n\varepsilon.
$$
\end{remark}

\section{Numerical applications}
\subsection{Linear steady problem}\label{sec:linear}
A first example illustrating the proposed approach based on $L^1$-norm minimization over dictionaries is considered here. It is the following parameterized advection equation with a source term:
\begin{equation}
\frac{\partial u}{ \partial t} - \frac{\partial u}{\partial x}(x) = f(x;\mu),~x\in[0,1],~t\in[0,T_{\max}], 
\end{equation}
where 
\begin{equation}
f(x;\mu) = \frac{2k\exp(-2k(x-\mu))}{(1+\exp(-2k(x-\mu)))^2}
\end{equation}
 $\mu\in[0.3,0.5]$ and $k=100$. 
Steady state solutions are considered next. For this configuration, the steady-states present a sharp gradient at location $x=\mu$. A Dirichlet boundary condition $u(0) = 1$ is imposed at $x=0$. The PDE is discretized by an upwind finite differences scheme using a uniform mesh, resulting in a HDM of dimension $N=10^3$.

\begin{figure}[h!]
\begin{center}
{\includegraphics[width=0.75\textwidth,clip=]{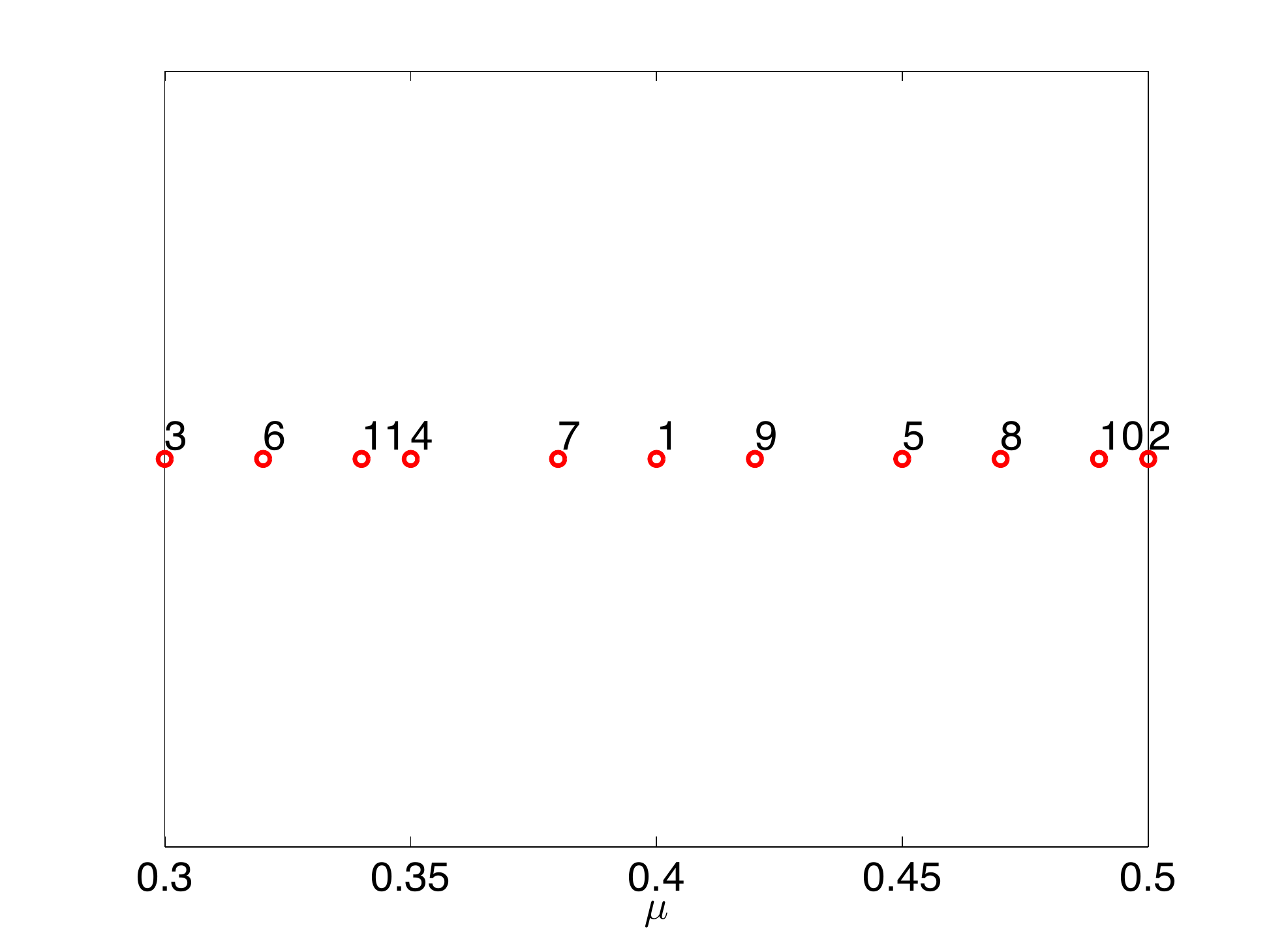}}
\end{center}
\caption{\label{fig:greedySamplesParam} Parameters selected by the greedy algorithm for the advection equation}
\end{figure}

\begin{figure}[h!]
\begin{center}
{\includegraphics[width=0.75\textwidth,clip=]{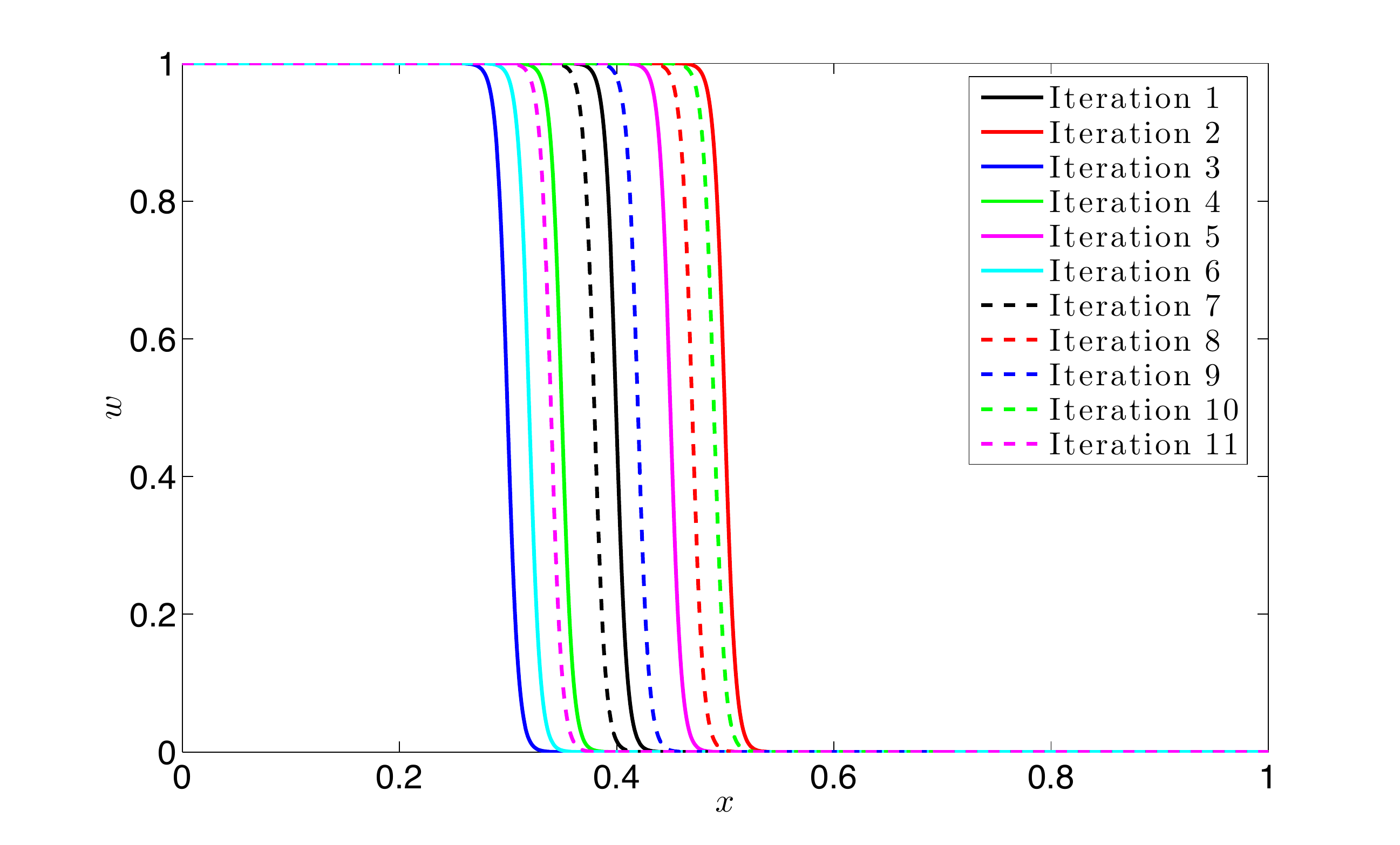}}
\end{center}
\caption{\label{fig:greedySamplesSnapshots} Dictionary members selected by the greedy algorithm for the advection equation}
\end{figure}

The greedy sampling algorithm proposed in Section~\ref{sec:training} is first applied to construct a dictionary that is accurate in the parametric domain $\mathcal{P} = [0.3,0.5]$. For that purpose, a set of $N_c=21$ candidate parameters $\mathcal{C}$ uniformly distributed in $\mathcal{P}$ is considered. An initial dictionary member $\mu_0=0.4$ is selected and then $10$ greedy iterations performed using a ROM based on Huber-norm minimization with the error indicator~\ref{eq:errorIndicL1}, resulting in a dictionary with $11$ members. The values of the parameters selected by the greedy approach are reported in Figure~\ref{fig:greedySamplesParam} and the dictionary members depicted in Figure~\ref{fig:greedySamplesSnapshots}. One can observe that the greedy algorithm selects in practice new samples that are maximally separated from the previously sampled dictionary members. 

\begin{figure}[h!]
\begin{center}
{\includegraphics[width=0.85\textwidth,clip=]{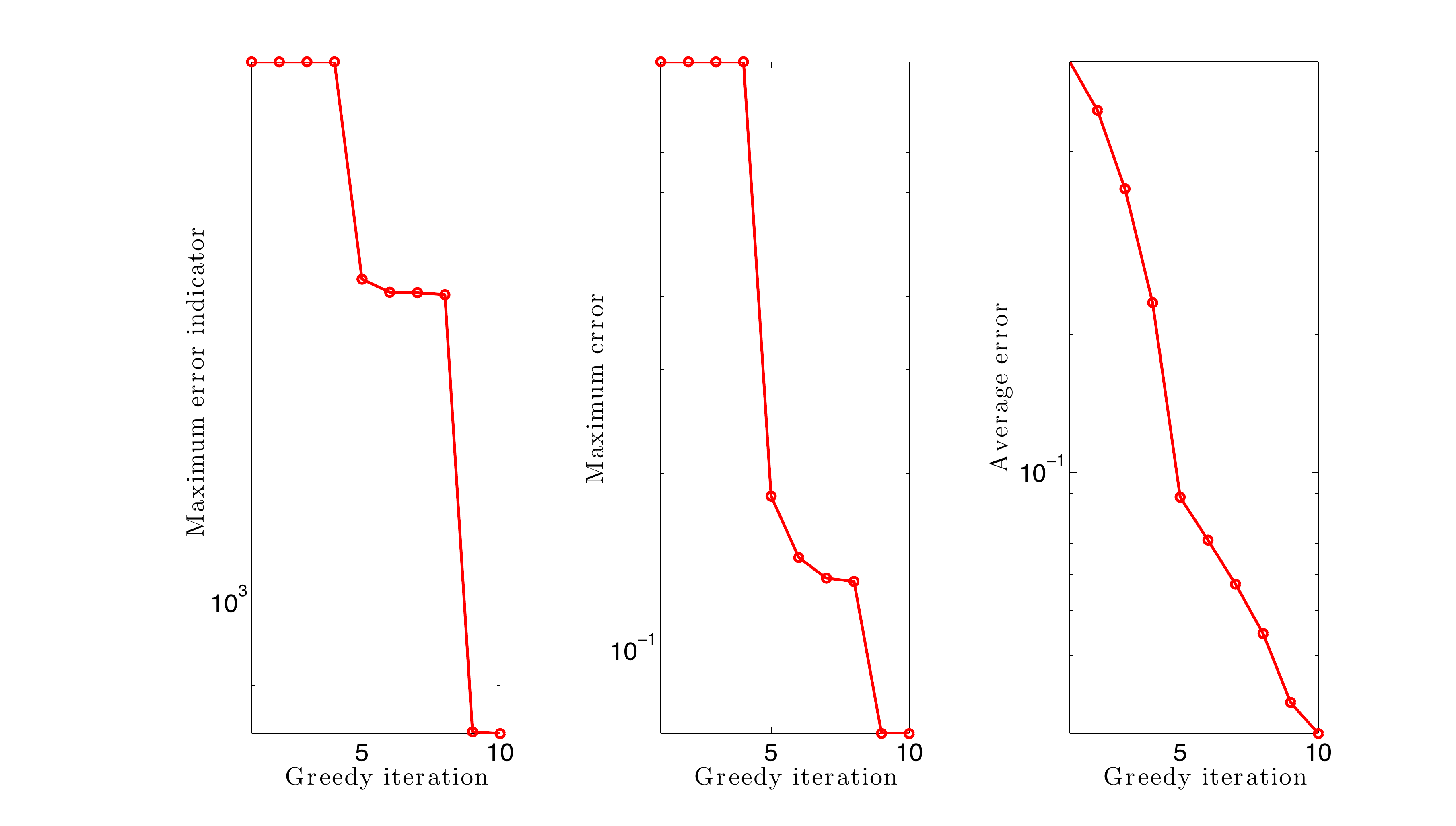}}
\end{center}
\caption{\label{fig:convergenceGreedy} Convergence of the greedy sampling procedure for the advection equation}
\end{figure}

\begin{figure}[h!]
\begin{center}
{\includegraphics[width=0.85\textwidth,clip=]{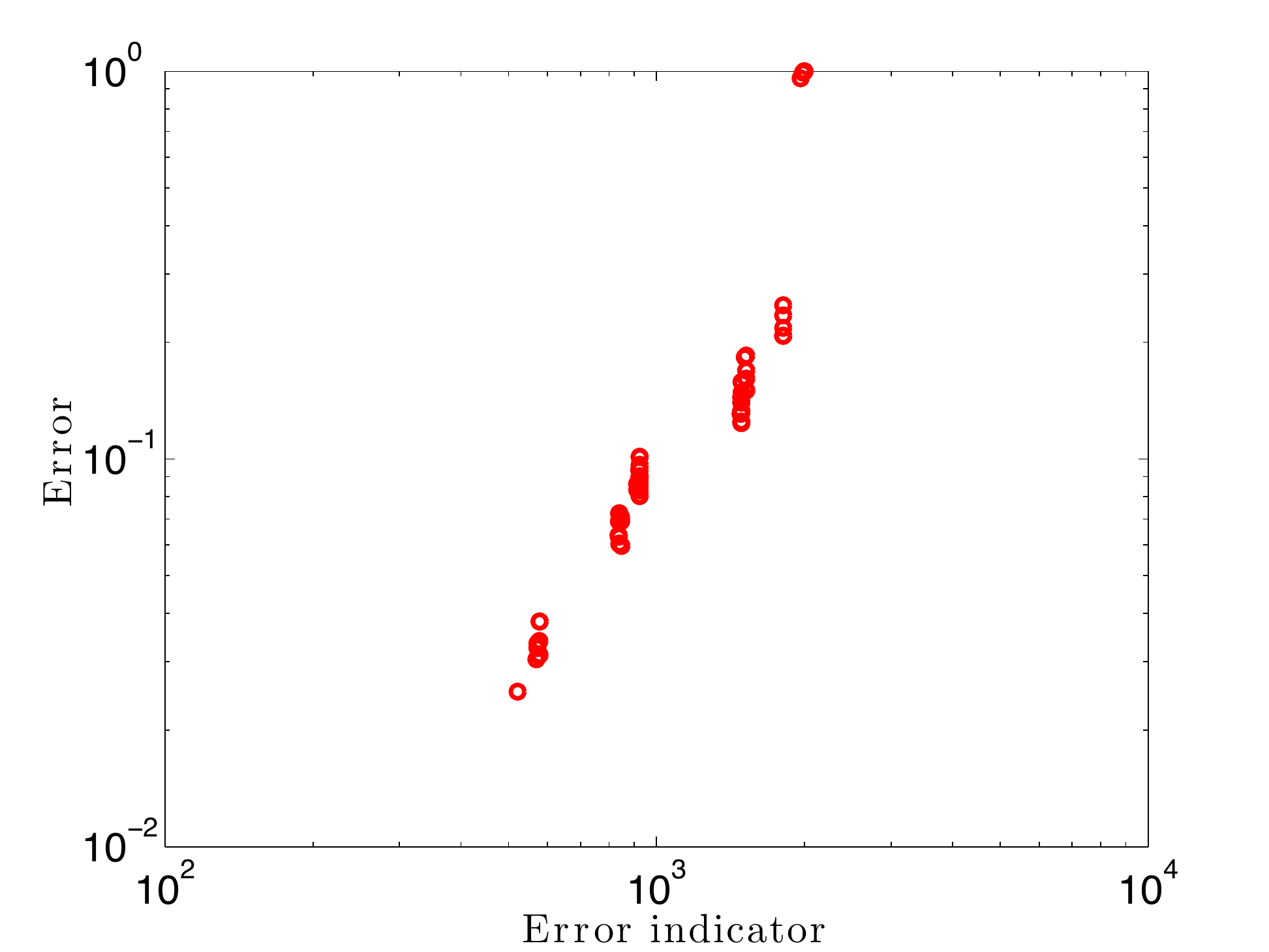}}
\end{center}
\caption{\label{fig:indicatorvserror} $L^2$-norm error as a function of the error indicator for the greedy iterates}
\end{figure}

A convergence analysis is presented in Figure~\ref{fig:convergenceGreedy} in which the maximum error indicator, maximum and average $L^2$-norm errors over the candidate sets are computed at each step of the greedy sampling procedure. One can observe that all errors are non-increasing with the iterations. Furthermore, the error and indicator follow a similar progression. This is confirmed by displaying the true error as a function of the indicator for all iterates of the greedy procedure in Figure~\ref{fig:indicatorvserror}. One can observe a linear relationship between those two quantities that can be readily captured by a linear regression~\cite{pdt15}.

\begin{figure}[h!]
\begin{center}
{\includegraphics[width=0.85\textwidth,clip=]{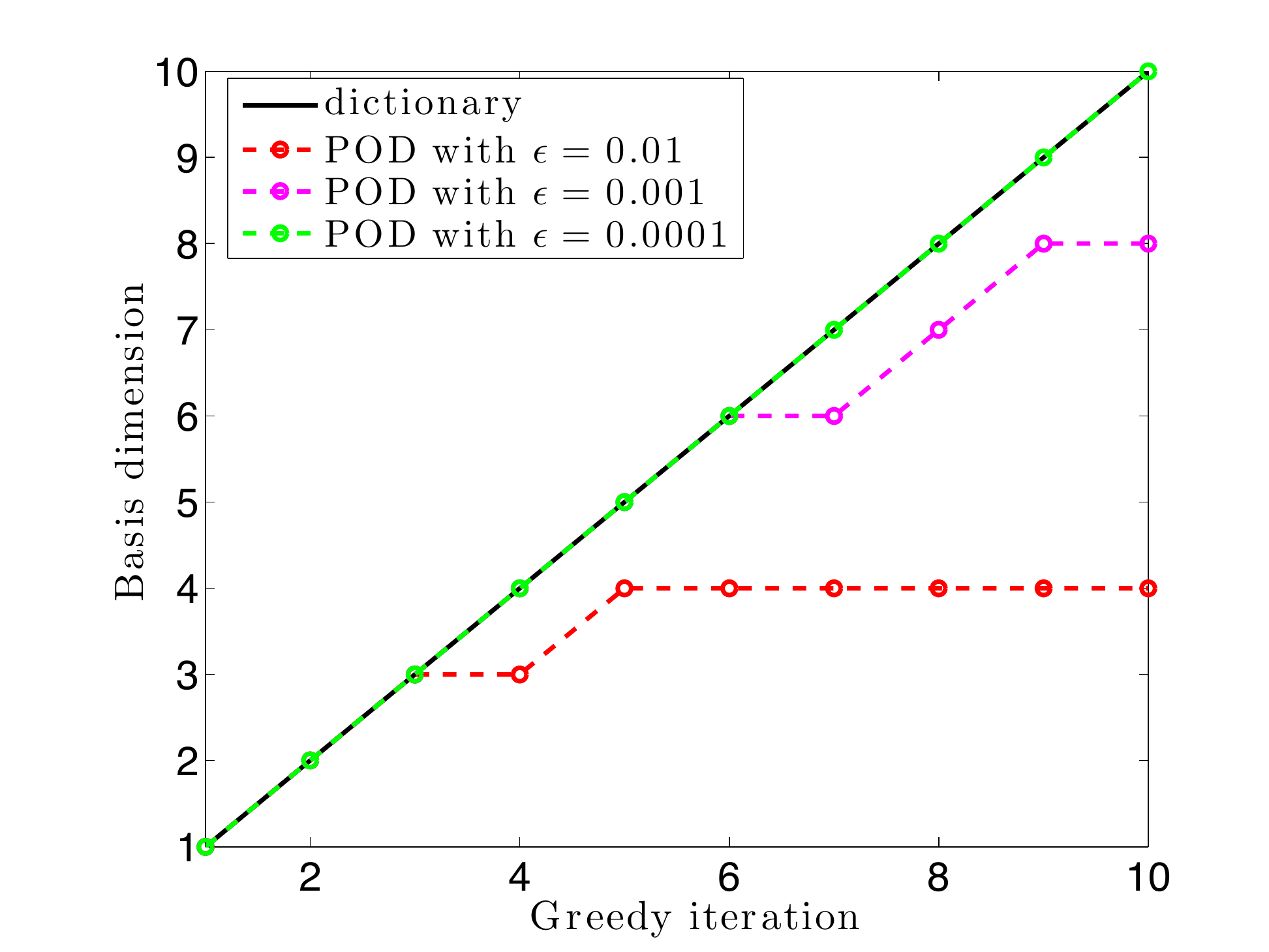}}
\end{center}
\caption{\label{fig:ROBdim} Dimension of the reduced bases constructed by the dictionary procedure and POD}
\end{figure}

\begin{figure}[h!]
\begin{center}
{\includegraphics[width=0.85\textwidth,clip=]{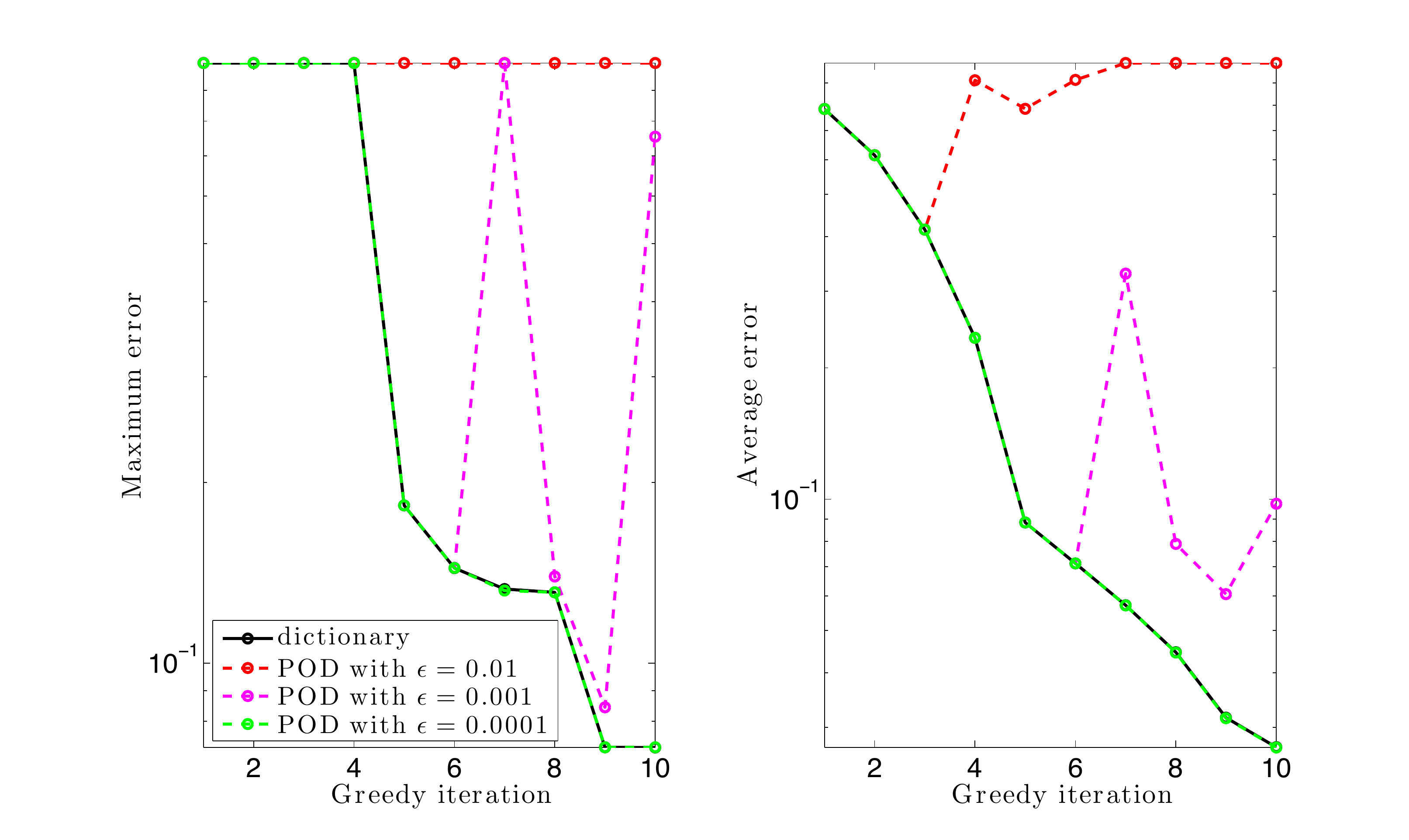}}
\end{center}
\caption{\label{fig:dicovsPOD} Errors associated with the dictionary procedure and POD}
\end{figure}

Next, the dictionary approach ---in which all the content of the sampled snapshots is retained--- is compared to POD for which the information associated with the singular vector with smallest energy is truncated. For that purpose, POD is applied with three levels of energy truncation $\epsilon_{\text{POD}} \in \{10^{-2},10^{-3},10^{-4}\}$ to the sets of snapshots obtained at each iteration of the greedy sampling procedure.  One can first observe from Figure~\ref{fig:ROBdim} that the energy levels $\epsilon_{\text{POD}}\in\{10^{-2},10^{-3}\}$ result in truncation of the content of the snapshots, while applying POD with $\epsilon_{\text{POD}}=10^{-4}$ does not result in truncation and is equivalent to the dictionary procedure. The maximum and average error over the candidate sets are displayed in Figure~\ref{fig:dicovsPOD}. Furthermore, when truncation occurs, POD results in much larger errors as well as several cases in which adding a snapshot to the sampled set increases the error. On the other hand, the dictionary approach results in a monotonous decrease of error.

\begin{figure}[h!]
\begin{center}
{\includegraphics[width=0.85\textwidth,clip=]{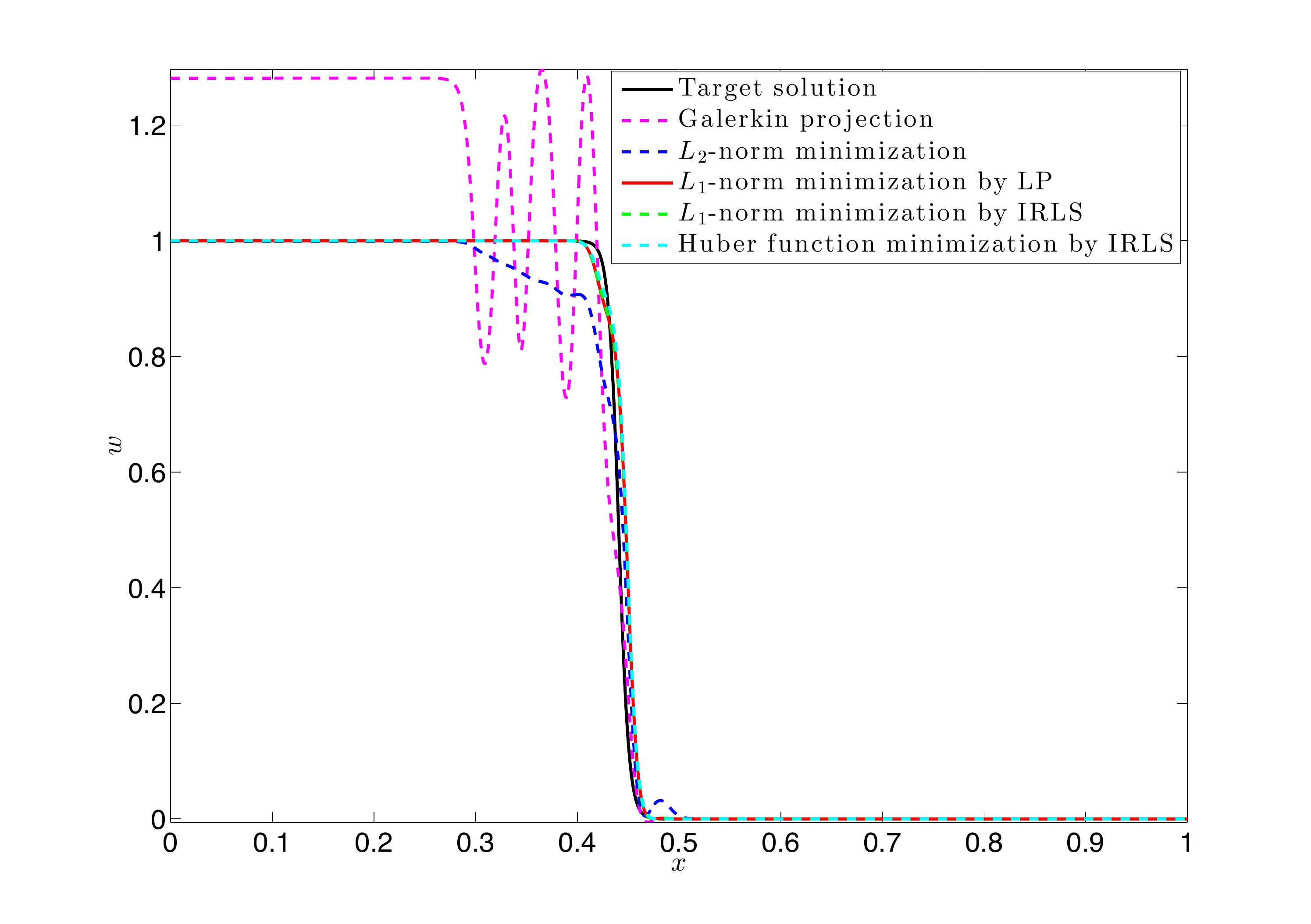}}
\end{center}
\caption{\label{fig:allMOR} Comparison of all MOR approaches and the target solution at $\mu^\star = 0.4412$}
\end{figure}

Finally, a target parameter $\mu^\star = 0.4412$ is randomly selected and the dictionary approach based on the previously constructed $11$ sampled is tested together with the following five model reduction approaches:
\begin{enumerate}
\item Galerkin projection
\item Minimization of the $L^2$-norm of the residual 
\item Minimization of the $L^1$-norm of the residual by Linear Programming
\item Minimization of the $L^1$-norm of the residual by IRLS with tolerance $\epsilon=10^{-4}$
\item Minimization of the Huber function applied to the residual by IRLS with tolerance $\epsilon=10^{-4}$
\end{enumerate}
The solutions obtained using each MOR approach are compared to the target solution in Figure~\ref{fig:allMOR}. Qualitatively, one can observe that the $L^1$-norm and Huber function-based approaches approximate the target solution the best by providing solutions with steep discontinuities. On the other hand, Galerkin projection leads to a very oscillating solution that is completely inaccurate. Minimization of the $L^2$-norm of the residual leads to a solution that presents undershoot and overshoots before and after the discontinuity, respectively. The relative $L^2$-norm errors between the target solution and each ROM solution are reported in Table~\ref{tab:Errors}. One can observe that the approaches based on $L^1$-norm minimization (including the Huber function) lead to the smallest errors. In that same table, the CPU timings are reported. One can observe that the IRLS procedure is less computationally expensive that the Linear Programing approach by more than one order of magnitude. The Huber function minimization approach is slightly more expensive than the $L^2$-norm minimization and Galerkin approaches, but leads to much more accurate predictions, as observed in Figure~\ref{fig:allMOR}. It is also less computationally expensive than the $L^1$-norm minimization approach. This is due to the fact that the IRLS procedure required $43$ iterations to converge in the case of $L^1$-norm minimization while for the Huber function it only required $12$ iterations to converge to the same tolerance. This discrepancy may be attributed to the non-differentiability of the $L^1$-norm at zero which may leads to a slower convergence. 

\begin{figure}[h!]
\begin{center}
{\includegraphics[width=0.85\textwidth,clip=]{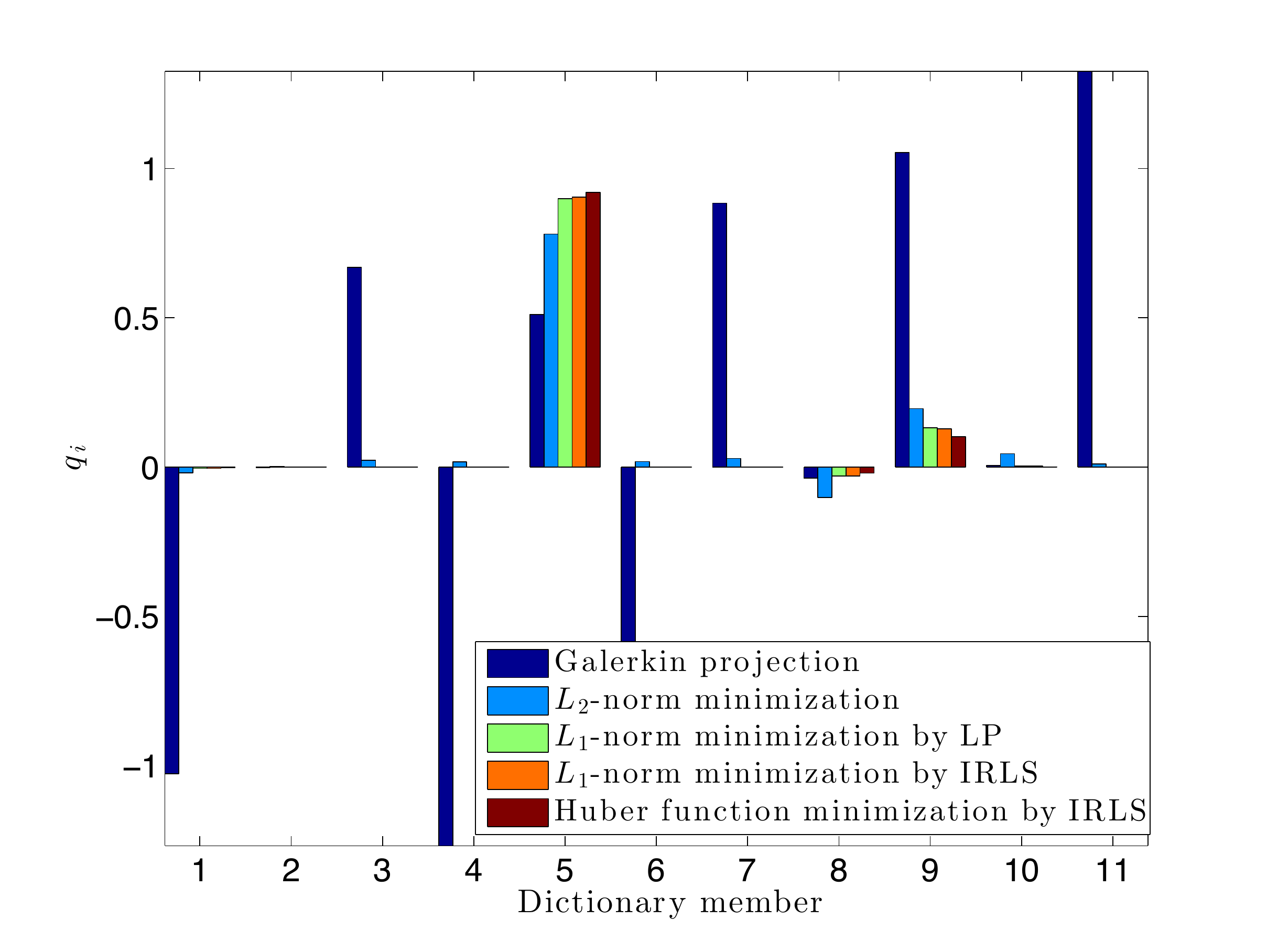}}
\end{center}
\caption{\label{fig:redCoord} Reduced coordinates of the solutions for all MOR approaches at $\mu^\star = 0.4412$}
\end{figure}

Finally, the reduced coordinates associated with each ROM are reported in Figure~\ref{fig:redCoord}. The $L^1$-norm and Huber function minimizations lead to sparse solutions whether Galerkin projection and $L^2$-norm minimization have all non-zero contributions from all dictionary members.

\begin{table}[h!]
\begin{center}\begin{tabular}{|c|c|c|c|c|c|}
\hline
 & Galerkin  & $L^2$-norm  & $L^1$-norm (LP) & $L^1$-norm (IRLS) & Huber function (IRLS)  \\
 \hline
 Relative error &0.2540 & 0.0600 & 0.0544 & 0.0549 & 0.0584 \\
 \hline
  CPU timings (s) &0.0112 & 0.0101 & 0.7826 & 0.0681 & 0.0132 \\
 \hline
 \end{tabular} 
 \caption{One-dimensional advection equation: relative errors}
 \label{tab:Errors}
\end{center}
\end{table}

\subsection{Nonlinear steady problem}
\subsubsection{Unsteady Burgers' equation}\label{Unsteady Burgers}\label{sec:burgers}
\begin{figure}[h!]
\subfigure[Solutions]{\includegraphics[width=0.45\textwidth]{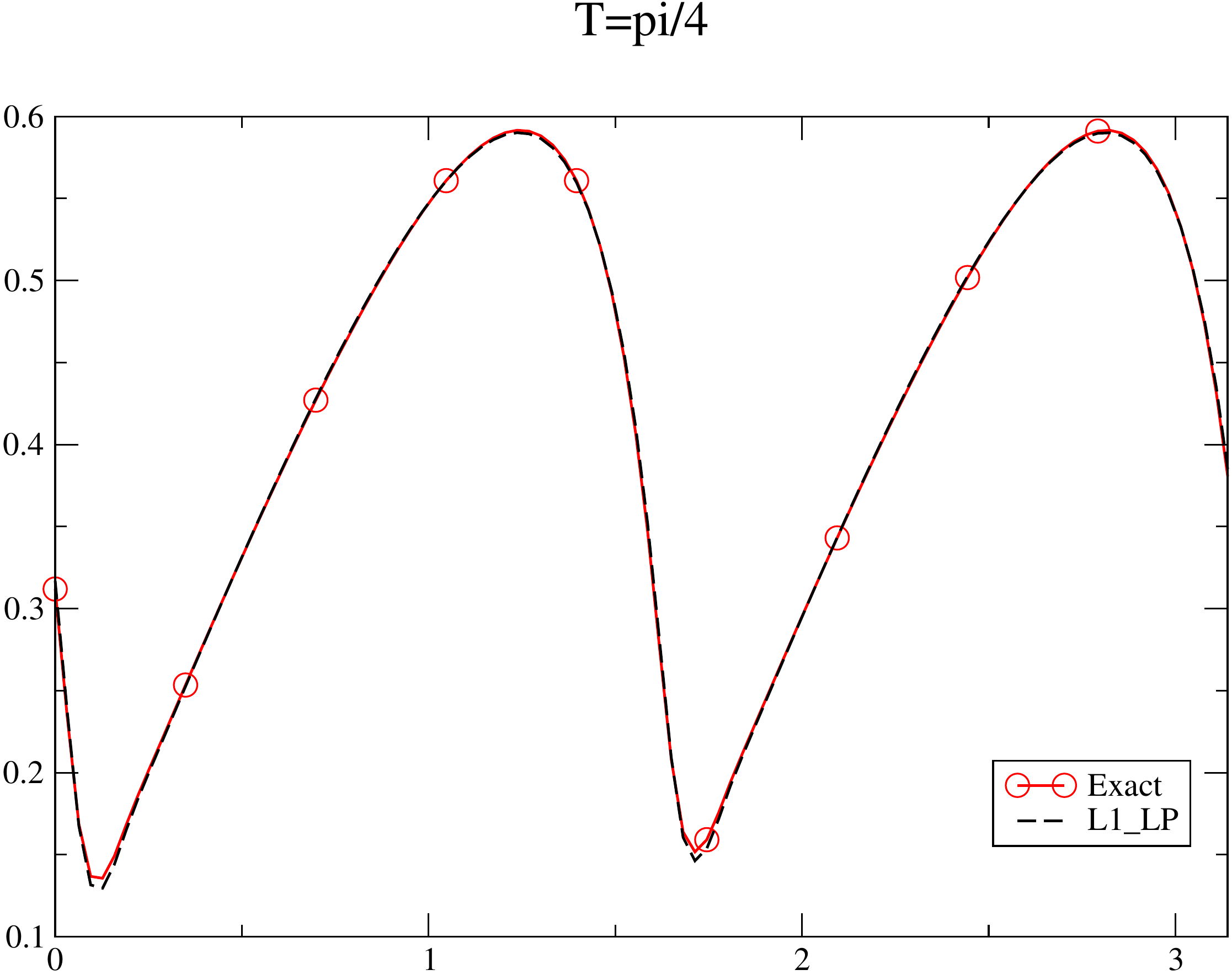}}
\subfigure[Zoom near a maximum]{\includegraphics[width=0.45\textwidth]{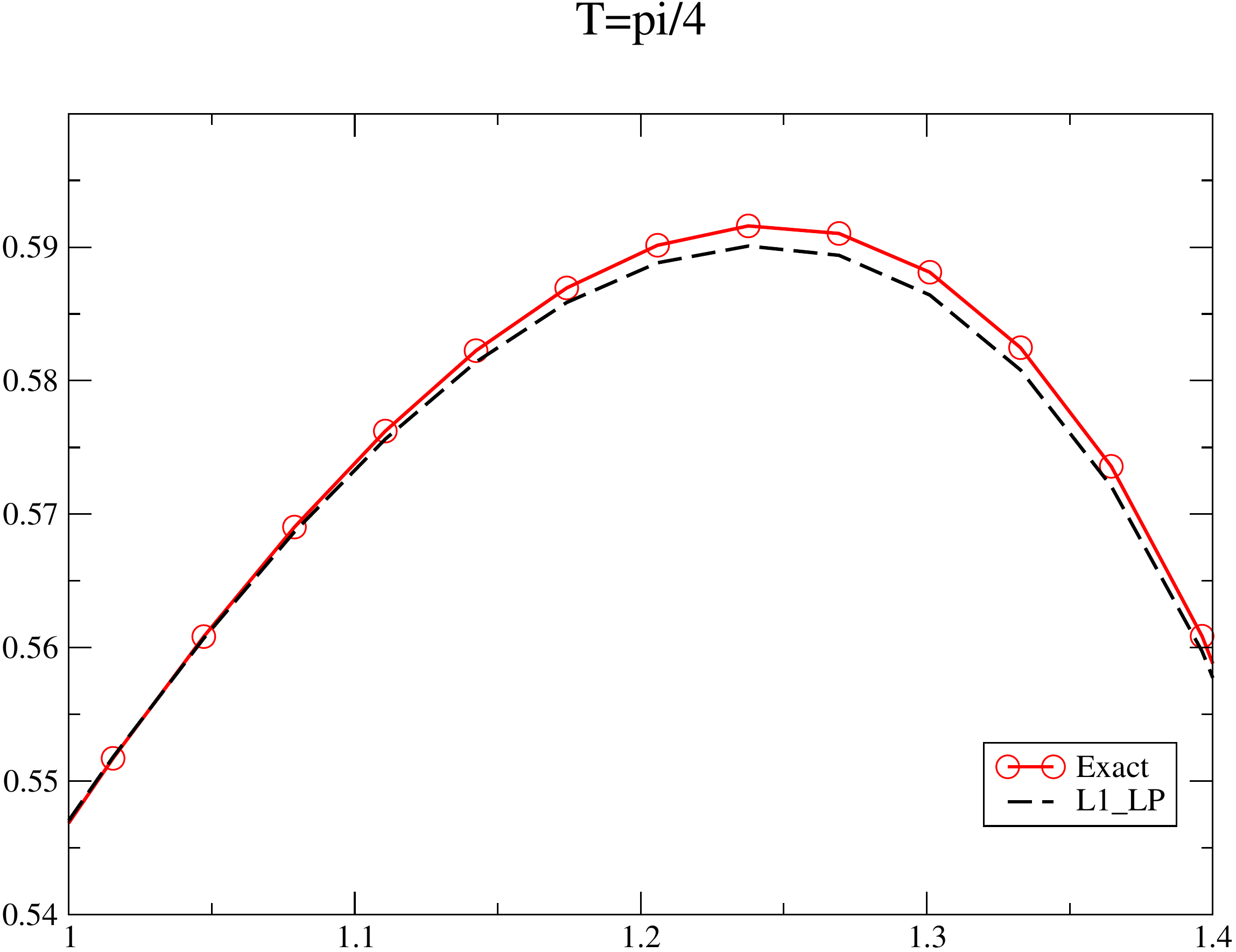}}
\centering{\subfigure[Zoom near a minimum]{\includegraphics[width=0.45\textwidth]{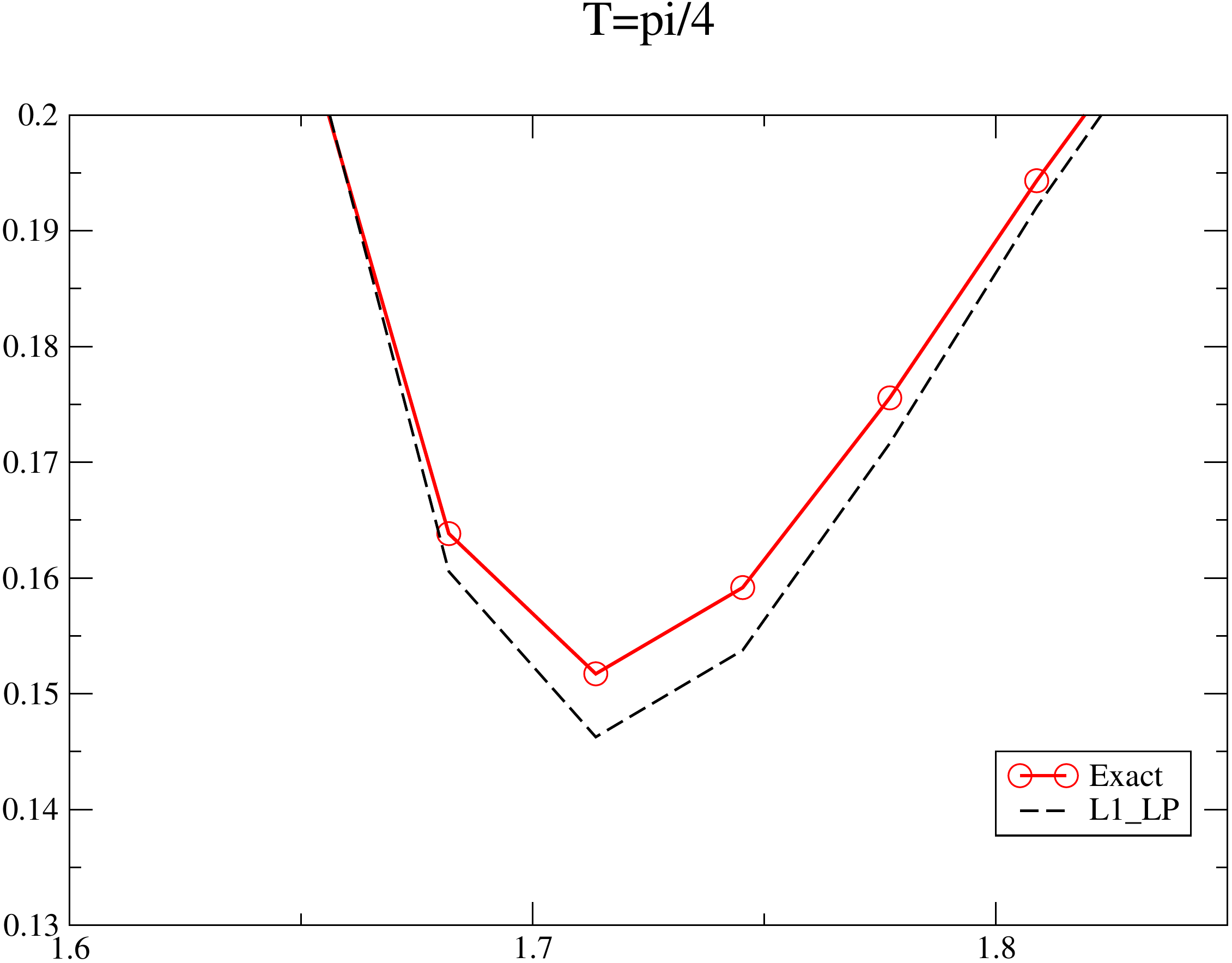}}}
\caption{\label{burger:pis4} Unsteady Burgers' equation: predicted solutions at target parameter $\mu^\star=0.5$ at $t=\frac{\pi}{4}$ }
\end{figure}

\begin{figure}[ht]
\subfigure[]{\includegraphics[width=0.45\textwidth]{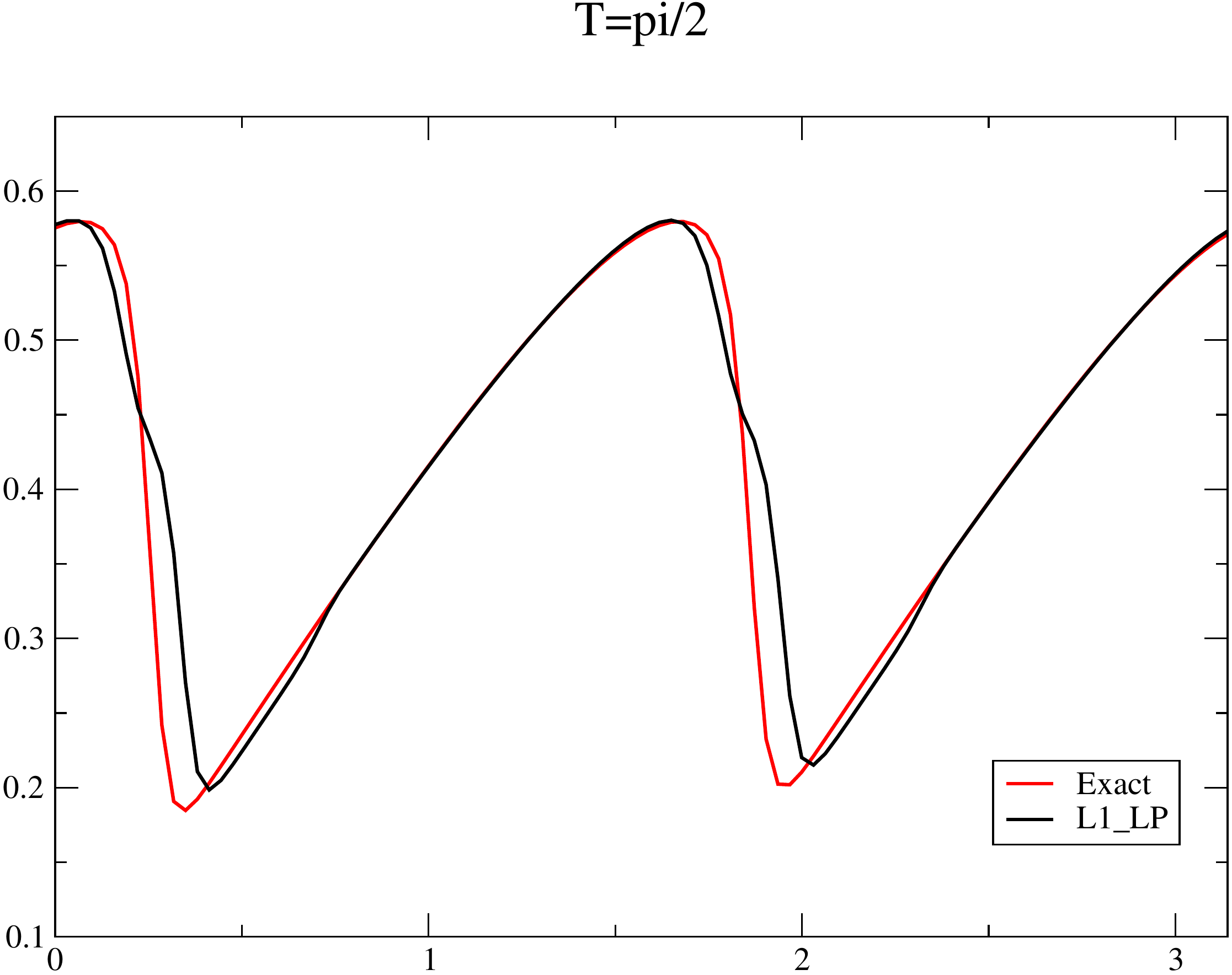}}
\subfigure[]{\includegraphics[width=0.45\textwidth]{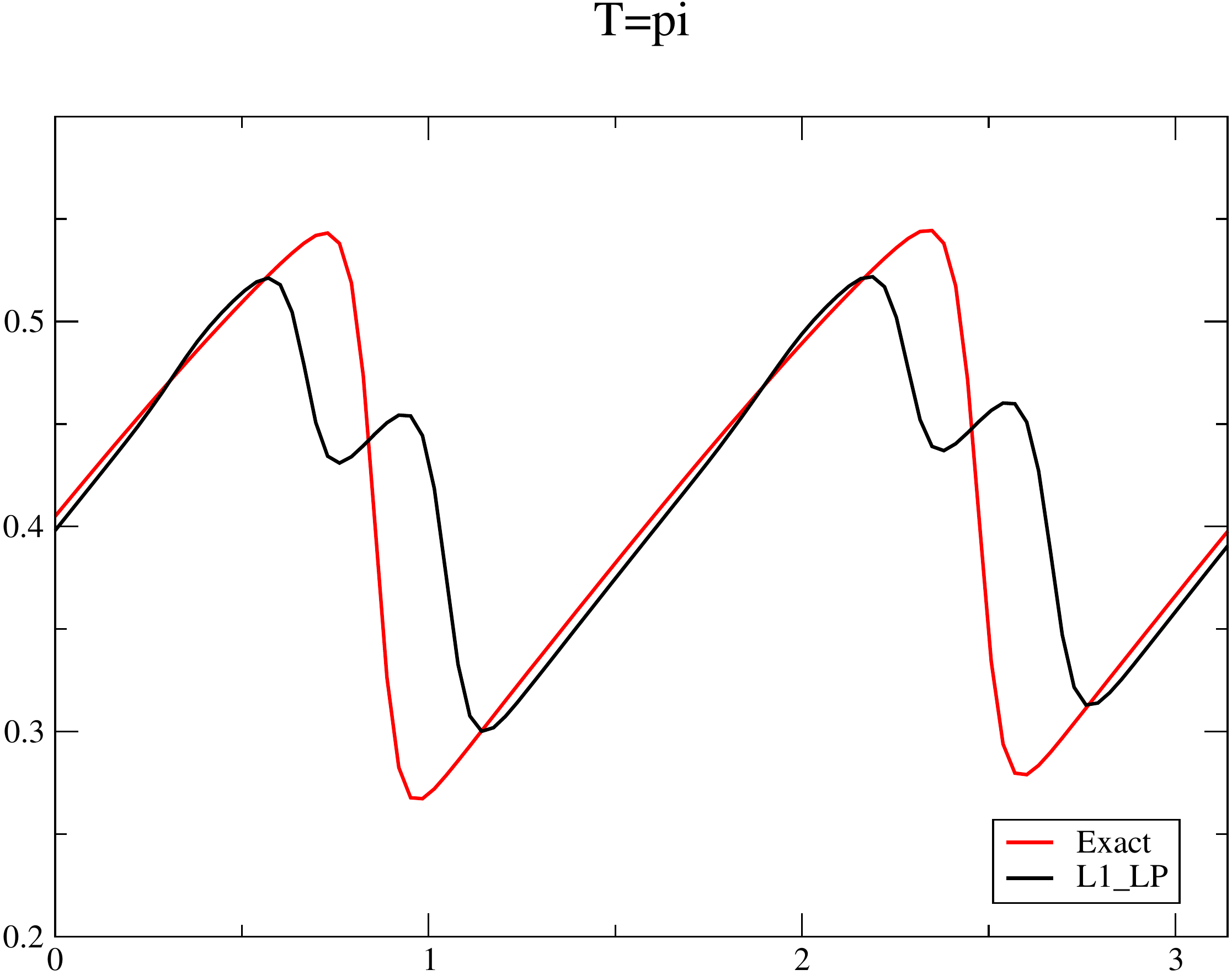}}
\caption{\label{burger:apres_choc} Unsteady Burgers' equation: predicted solutions at target parameter $\mu^\star=0.5$ at $t=\frac{\pi}{2}$ (left) and $t=\pi$ (right)}
\end{figure}

We consider here the system~\eqref{eq:2} in $\Omega=[0,2\pi]$ with periodic boundary conditions and the initial conditions parameterized  by
$$
u_0(x;\mu)=\mu \; \big | \sin(2\; x)\big |+0.1,
$$
where $\mu \in [0,1]$. In this setting, the solution develops a shock that moves with the velocity $\sigma_\mu=0.6 \mu$.
A dictionary $\mathcal{D}$ is constructed by sampling the parameters $\{0,0.2,0.4,0.6,1.0\}$ ($r=5$) and the solution sought for the predictive case $\mu^\star=0.5$. A shock appears at $t=1$. We display the solutions obtained by $L^1$-norm by LP minimization procedure for $t=\tfrac{\pi}{4}<1$, $t=\tfrac{\pi}{2}$ and $t=\pi$ in Figures \ref{burger:pis4} and \ref{burger:apres_choc}.

After the shock, the $L^1$-norm-type solutions are all close to each other and the shock is rather well reproduced with, however, an artifact that develops for longer  times, as seen at $t=\pi$. Nevertheless, the $L^1$-norm-type solutions are within the bounds of the ``exact" solution, and no large oscillation develops.

In a second set of numerical experiments, we consider the influence of the sampling parameter set included in the dictionary $\mathcal{D}$. We consider two dictionaries
 $\mathcal{D}_1=\{0.4,0.45, 0.55$ $ , 0.6\}$ and  $\mathcal{D}_0=\{0,0.2,0.4,0.45, 0.55, $ $ 0.6,1.0\}$, for the same target value of $\mu^\star=0.5$. These choices amounts to selecting samples close to the target value $0.5$ while varying elements of the dictionary that are not close to $0.5$.
 
\begin{figure}[h!]
\subfigure[$\mathcal{D}_1$]{\includegraphics[width=0.45\textwidth]{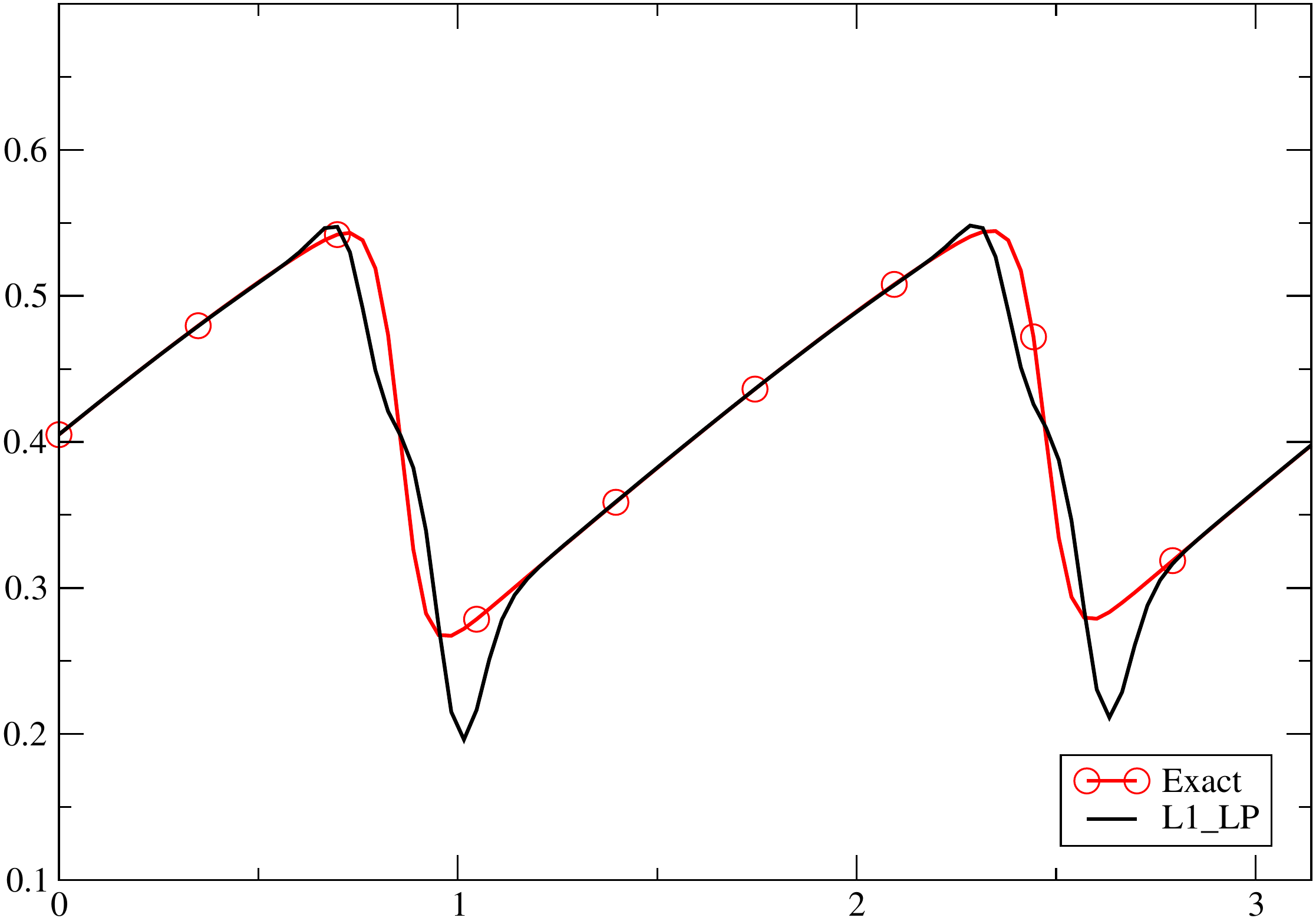}}
\subfigure[$\mathcal{D}_0$]{\includegraphics[width=0.45\textwidth]{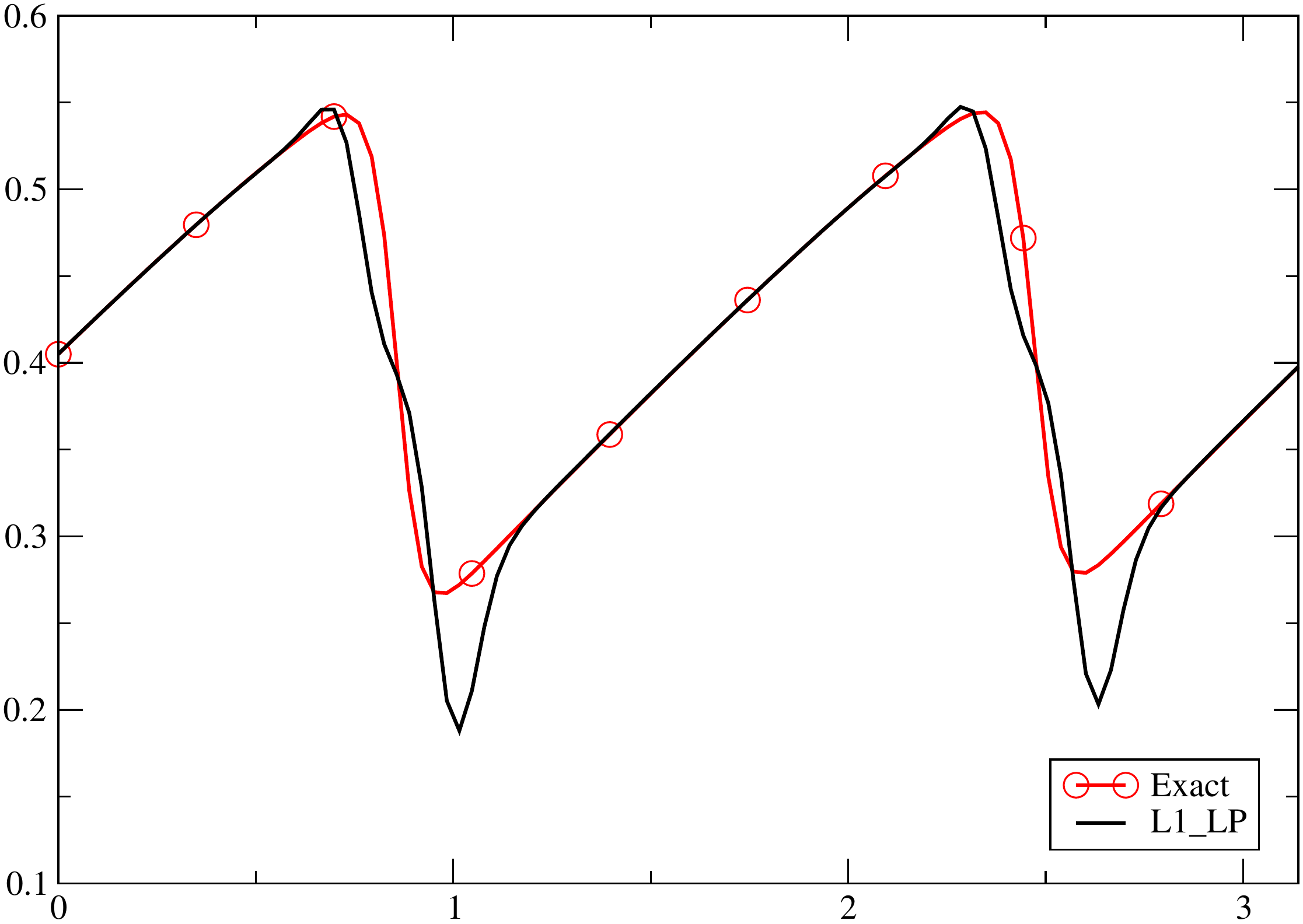}}
\caption{\label{influence} Unsteady Burgers' equation: predicted solutions at target parameter $\mu^\star=0.5$ at $t=\pi$ for two dictionaries associated with two samples of the parameter domain $\mathcal{P}$}
\end{figure}

We see that refining the dictionary has a positive influence as the target solution is much closer to the dictionary elements. This is confirmed by additional experiments where the samples of  $\mu$ used to generate the dictionary where more numerous and closer to $0.5$ (not reported here). The $L^1$-norm-type solutions are however  unaffected by the presence of these ``outliers" in the dictionary. 

\subsubsection{Euler equations}\label{sec:euler}
The one-dimensional Euler equations are considered on $\Omega = [0,1]$
\begin{subequations}\label{euler}
\begin{equation}\label{euler:1}
\dpar{}{t}\begin{pmatrix}\rho \\ \rho u \\E\end{pmatrix}+\dpar{}{x}\begin{pmatrix}
\rho u \\ \rho u^2+p \\ u(E+p)\end{pmatrix}=0,
\end{equation}
for which $U=(\rho,\rho u,E)^T$ and the pressure is given by
\begin{equation}
\label{euler:2}
p=(\gamma-1)\left( E-\frac{1}{2}\rho u^2\right)
\end{equation}
with $\gamma=1.4$.

This problem is parametrized by the initial conditions $U_0(x;\mu)$. To define the parametrized initial conditions of the problem, the Lax and Sod cases are first introduced as follows. 

The state  $U_{\text{Sod}}(x)$ is defined by the primal physical quantities:
\begin{equation}
\label{euler:4}
V_{\text{Sod}}(x) =\left\{\begin{array}{ll}
\rho= 1 \text{ if }x\leq 0.5, & 0.125 \text{ otherwise,}\\
u=0.0 \\
p=1.0 \text{ if }x\leq 0.5, & 0.1 \text{ otherwise,}
\end{array}\right .
\end{equation}
and $U_{\text{Lax}}(x)$ defined by
\begin{equation}
\label{euler:5}
V_{\text{Lax}}(x)=\left\{\begin{array}{ll}
\rho= 0.445 \text{ if }x\leq 0.5, & 0.5 \text{ otherwise,}\\
u=0.698 \text{ if }x\leq 0.5, & 0.0 \text{ otherwise,}\\
p=3.528 \text{ if }x\leq 0.5, & 0.571 \text{ otherwise.}
\end{array}\right .
\end{equation}

The Sod condition presents a fan, followed by a contact and a shock. For the density and the pressure, the solution behaves monotonically, and the contact is moderate. The Lax solution has a very different behavior and the contact is much stronger. This is depicted in Figure \ref{figure:sod-lax} where the two solutions are shown for $t=0.16$.
\begin{figure}[h!]
\begin{tabular}{cc}
\includegraphics[width=0.45\textwidth]{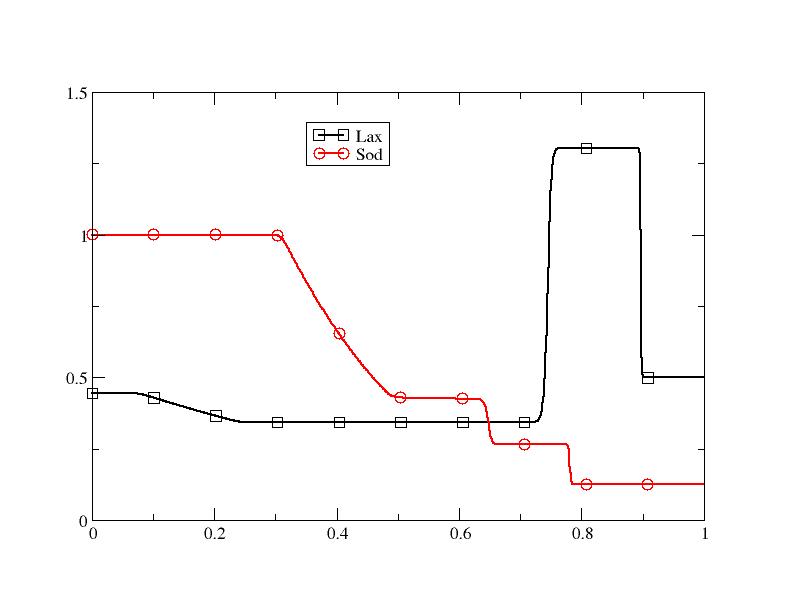}&  \includegraphics[width=0.45\textwidth]{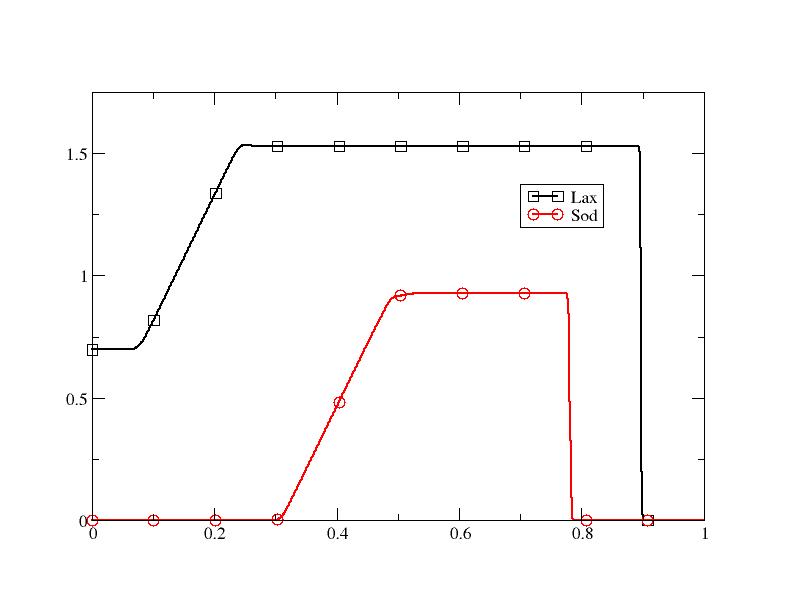}\\
Density& Velocity\\
\multicolumn{2}{c}{\includegraphics[width=0.45\textwidth]{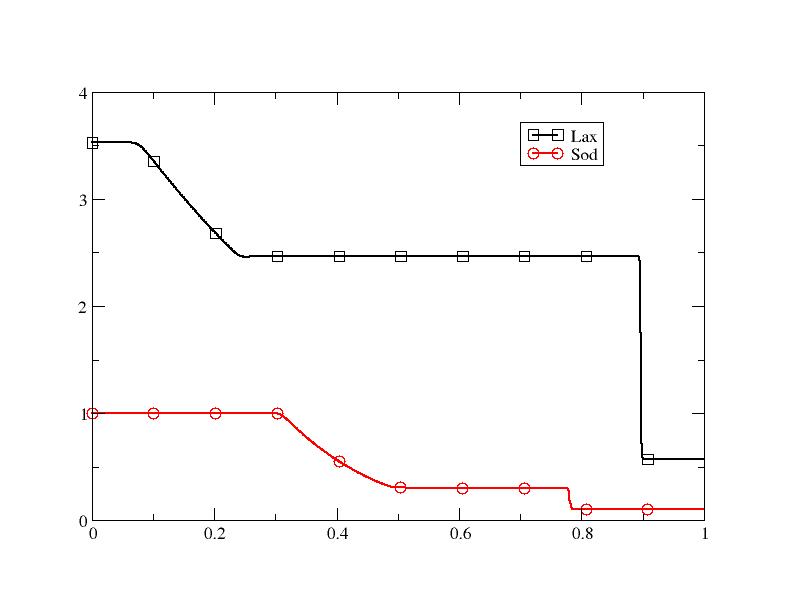}}\\
\multicolumn{2}{c}{Pressure}
\end{tabular}
\caption{\label{figure:sod-lax} One-dimensional Euler equations: density, velocity and pressure for the Lax and Sod problems}
\end{figure}

The initial condition are parametrized for $\mu\in [0,1]$ as 
\begin{equation}
\label{euler:3}
V_0(x;\mu)=\mu V_{\text{Sod}}(x)+(1-\mu)V_{\text{Lax}}(x)
\end{equation}
\end{subequations}
and the conservative initial variables $U_0(x;\mu)$ constructed from $V_0(x;\mu)$.

\begin{figure}[h!]
\subfigure[Density $\rho$]{\includegraphics[width=0.45\textwidth]{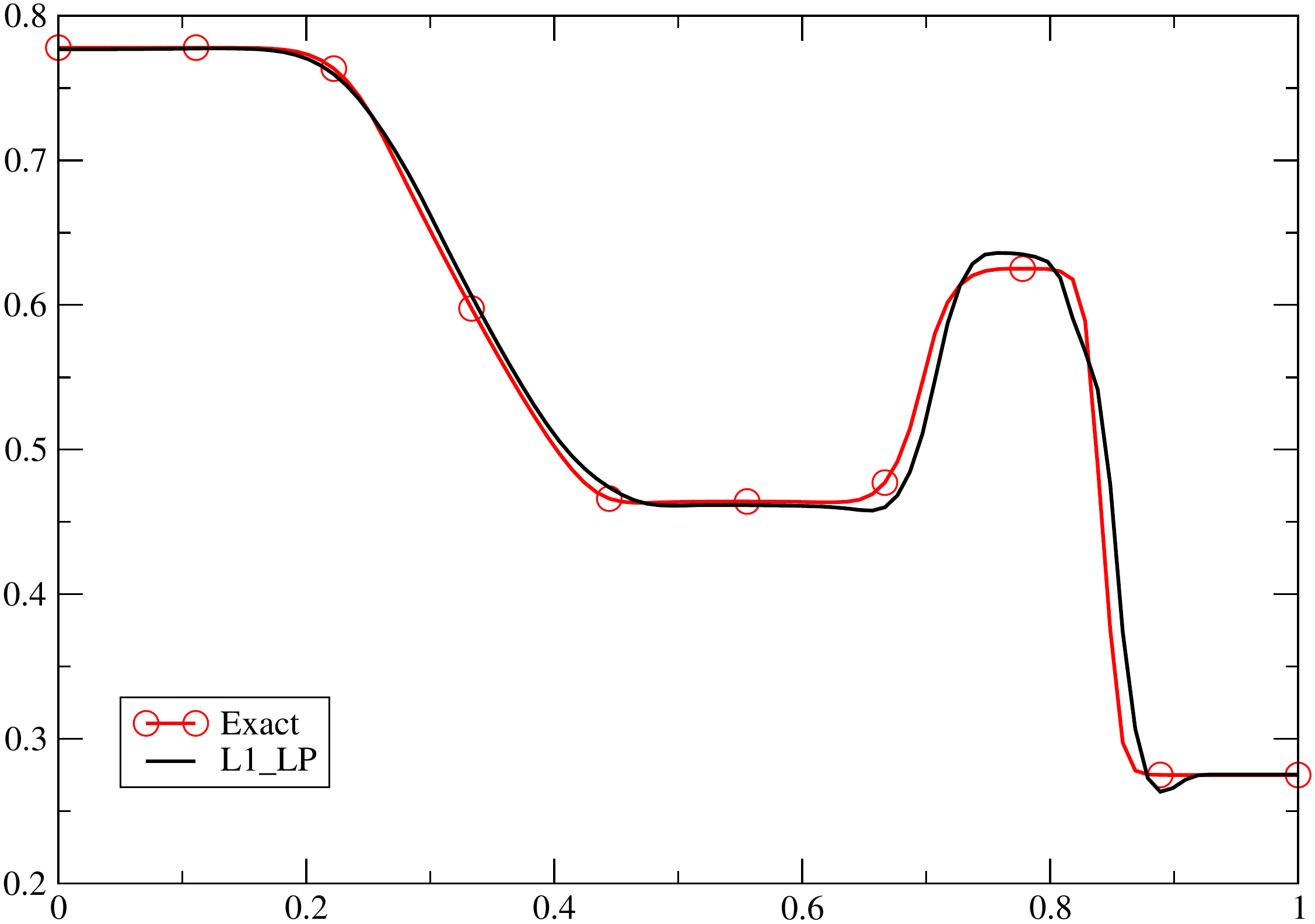}}
 \subfigure[Velocity $u$]{\includegraphics[width=0.45\textwidth]{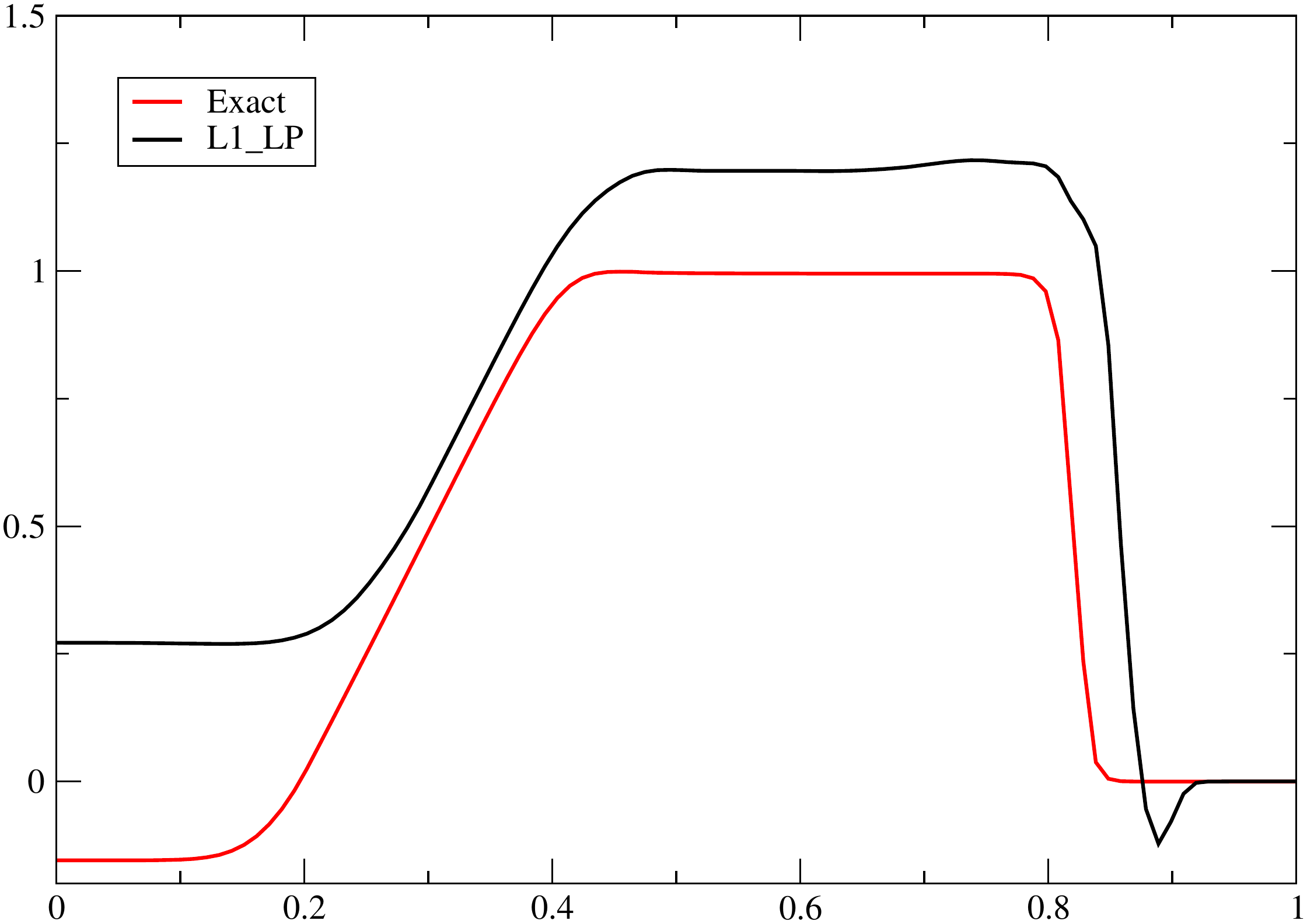} }
 \centering{\subfigure[Pressure $p$]{\includegraphics[width=0.45\textwidth]{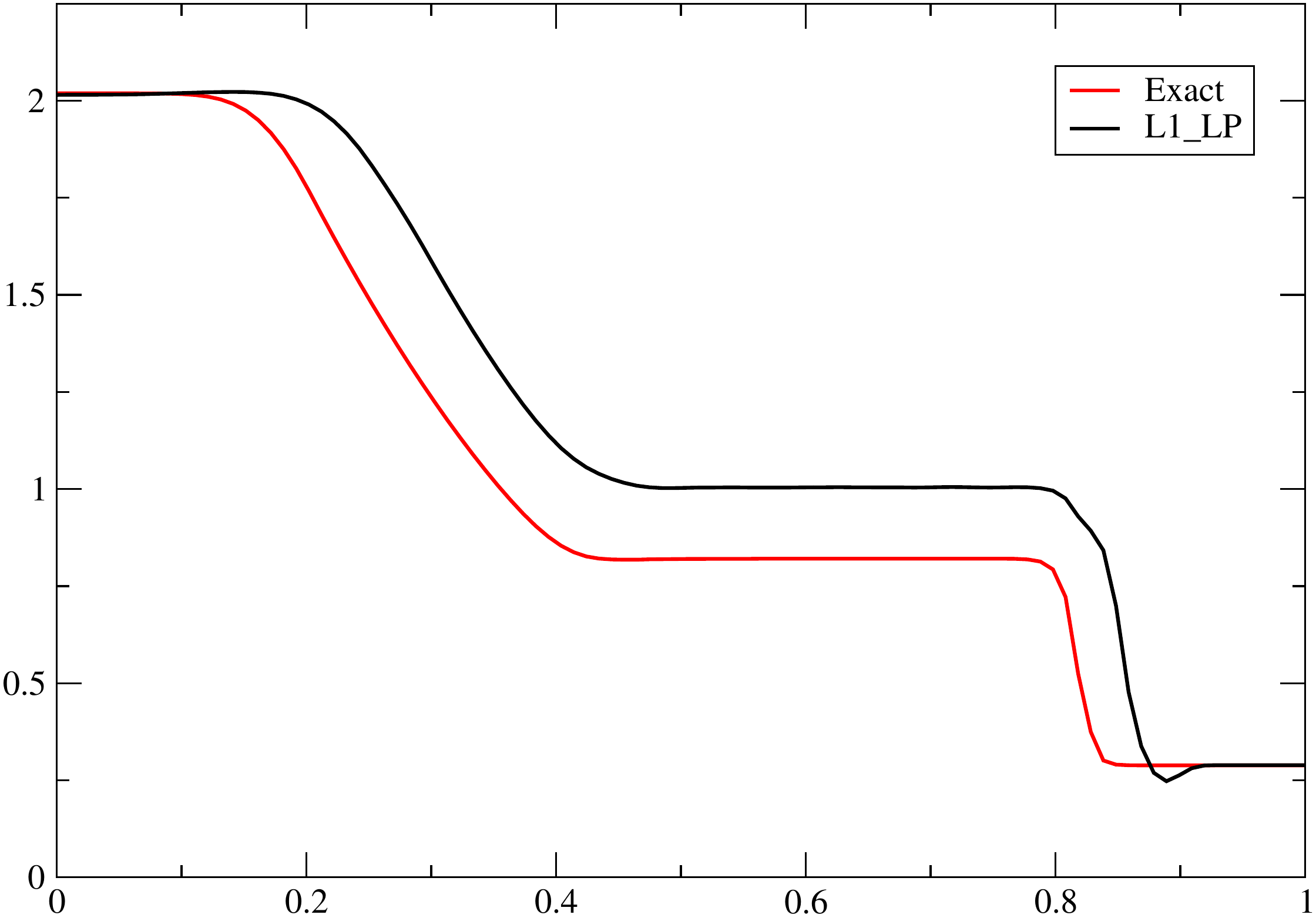} }}
 \caption{\label{fig:reconstruction:un} One-dimensional Euler equations: predicted solutions with strategy  \eqref{reconstruction:un} based on a single expansion}
 \end{figure}

\begin{figure}[h!]
\subfigure[Density $\rho$]{\includegraphics[width=0.45\textwidth]{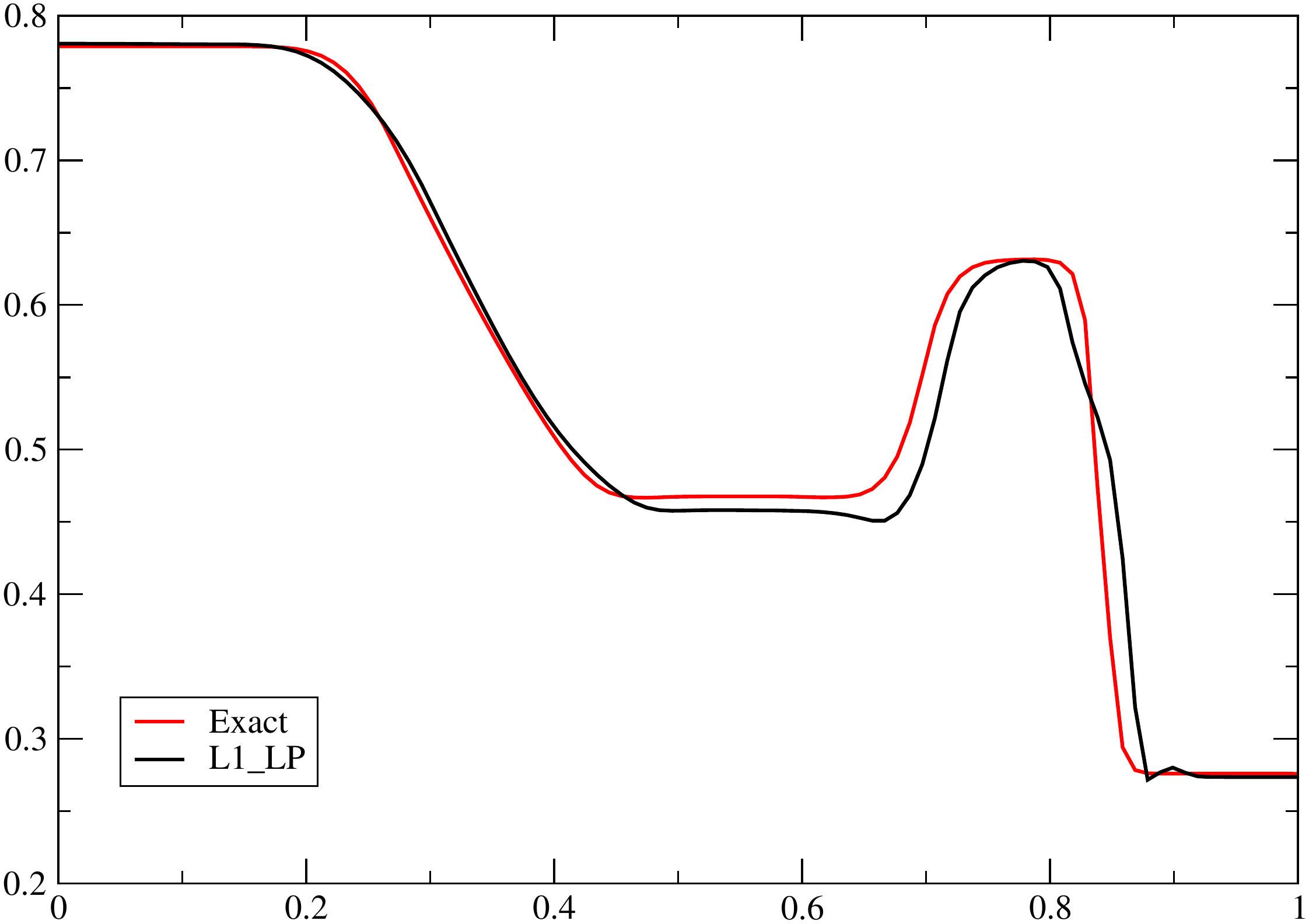}}
 \subfigure[Velocity $u$]{\includegraphics[width=0.45\textwidth]{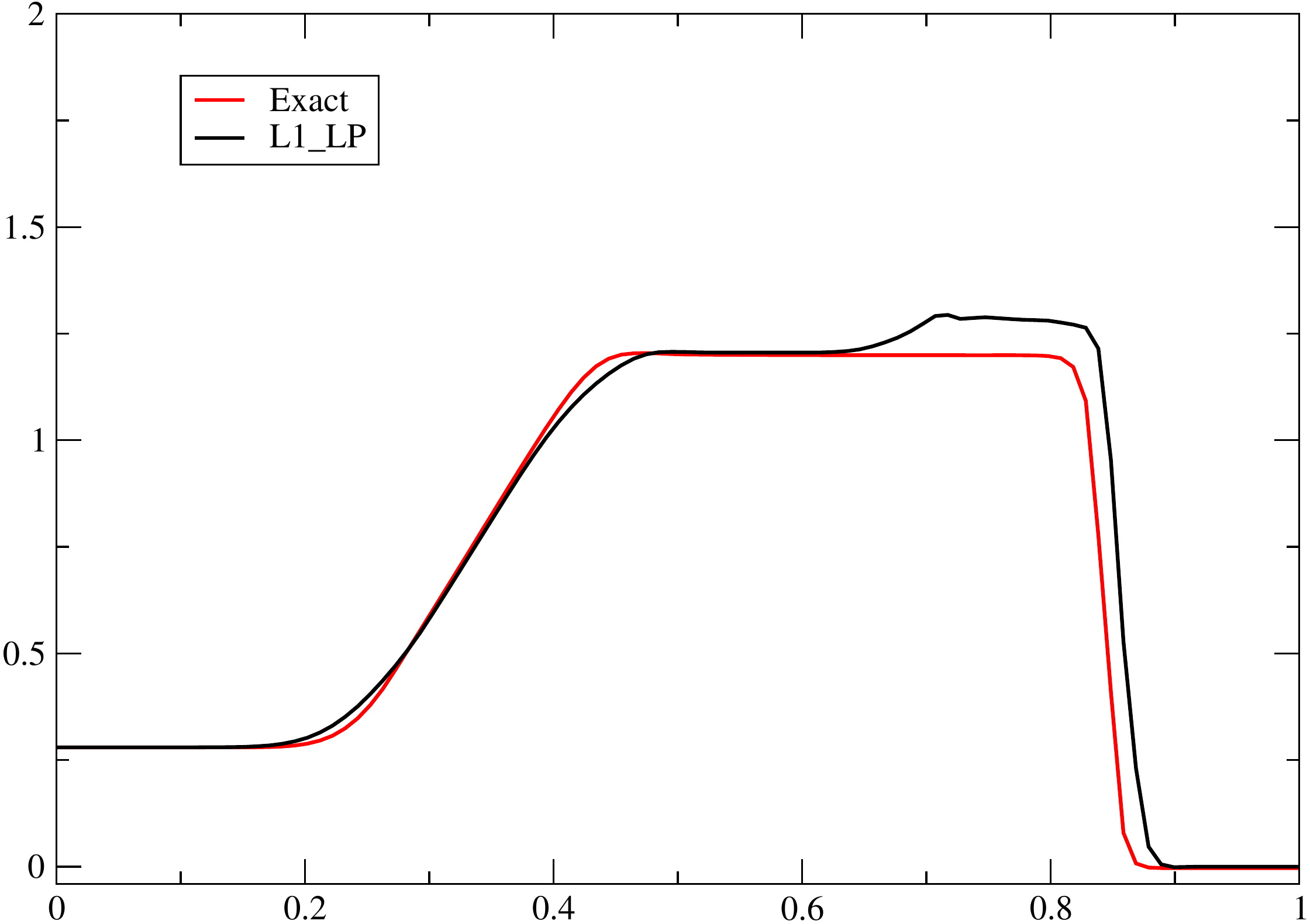}}
 \centering{\subfigure[Pressure $p$]{\includegraphics[width=0.45\textwidth]{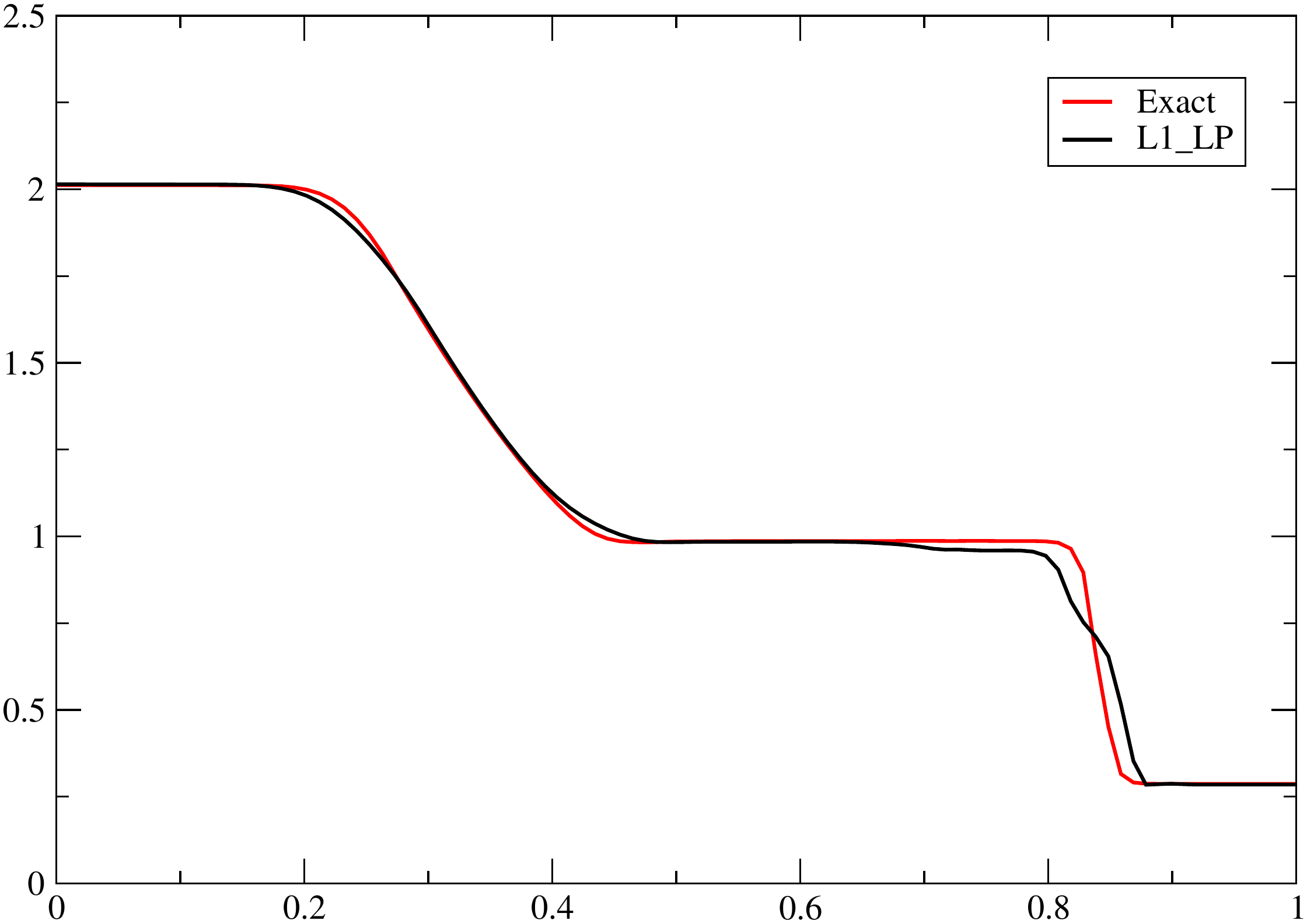}}}
 \caption{\label{fig:reconstruction:trois} One-dimensional Euler equations: predicted solutions with strategy \eqref{reconstruction:trois} based on multiple expansions}
\end{figure}

In the subsequent numerical experiments, two strategies are exploited to construct, from the dictionary $\mathcal{D}$, the approximation $\ubold^n(\mu)$ of the solution at each time step $n$:
\begin{itemize}
\item Either we reconstruct together the discretized density vectors $\rhobold$, momentum $\mbold=\rhobold \ubold$ and energy $\Ebold$, i.e. the state variable at time $t_n$ using only one coefficient vector $\alphabold^n=(\alpha^n_1,\cdots,\alpha^n_r)$
\begin{equation}
\label{reconstruction:un}
\ubold^n=\begin{pmatrix}
\rhobold^{n}\\
\mbold^n\\
\Ebold^n\end{pmatrix}\approx \sum_{j=1}^r \alpha_j^n  \ubold^n(\mu_j) . 
\end{equation}
Here the $\{\alpha_j^n\}_{j=1}^r$ are obtained by minimizing $J$  on the density components of the state because the density enable to detect fans, contact discontinuities and shocks, contrarily to pressure and velocity which are constant across contact waves. Doing so we expect to control better the numerical oscillations, if any, than with the other physical variables. Similar arguments could be applied with the other conserved variables as well.
\item Alternatively, we reconstruct each conserved variable separately
\begin{equation}\label{reconstruction:trois}
\rhobold^n\approx\sum_{j=1}^r \alpha_{j}^n \rhobold^n(\mu_j), \qquad
\mbold^n\approx\sum_{j=1}^r \alpha_{j}^n \mbold^n(\mu_j), \qquad
\Ebold^n\approx\sum_{j=1}^r \alpha_{j}^n \Ebold^n(\mu_j) .
\end{equation}
where the minimization procedures are done \textit{independently} on each conserved variable.
\end{itemize}
In order to test these approaches, the PDE is discretized by finite volumes using a discretization resulting in $Np=3000$ dofs. The parameter range $\mathcal{D}=\{0.0, 0.2,0.4,0.5,0.8,1\}$ is considered together with a target $\mu^\star=0.6$. The results using the first strategy, see eq. \eqref{reconstruction:un}, are displayed in Figure~\ref{fig:reconstruction:un} and those using the second strategy, see eq. \eqref{reconstruction:trois}, reported in Figure \ref{fig:reconstruction:trois}.
 
 From both figures, we can see that the overall structure of the solutions is correct. Nevertheless, there are differences that can be highlighted.
From Figure \ref{fig:reconstruction:un}, we can observe that the density predictions, besides an undershoot at the shock, are well reproduced. However, we cannot recover correct values of the initial velocity (see left  boundary), because there is no reason to believe that the coefficient $\alphabold$, evaluated from the density \emph{only}, will also be correct for the momentum. A careful observation of the pressure plot also reveals the same behavior which is not satisfactory. For the same reason, if any other \emph{single} variable is used for a global approximation of each conservative variables, there no reason why better qualitative results could be obtained.

This problem does not occur with the second strategy for the reconstruction \eqref{reconstruction:trois}: the correct initial values are recovered. We have some slight problems on the velocity, between the contact and the shock. 

In order to obtain these results we have been faced to the following issue. Take the momentum, for example. For at least half of the mesh points, its value is $0$, and for half of the points, its value is set to a constant. Hence, the matrix $\Abold$ used in the minimization procedure and built on the momentum dictionary has rank 2 only. The same is true for the other variables, and we are looking here for $r$ coefficients. Several approaches can be followed to address this issue. The first one relies on Gram-Schmidt orthogonalization of the solutions prior to their use as a basis for the solution. The second approach, followed here, consists into perturbing infinitesimally and randomly the matrices involved in the procedure,
so $A_{ij}$ is replaced by $A_{ij}+\varepsilon_{ij}$. The distribution of $\varepsilon_{ij}$ is uniform. This has the effect of giving the maximum possible rank to the perturbed matrix. We have expressed that $\epsilon_{ij}$ should depend on the variable, we have chosen
$$\varepsilon_{ij}=\epsilon_{ij} L_{\text{ref}}$$ where $L_{\text{ref}}$ is the difference between the minimum and the maximum, over the dictionary, of the considered variable. Choosing the same $\varepsilon_{ij}$ for all variables, this has the effect of increasing the amplitude of the oscillations after then shock.

All this being said, the solution using three distinct coefficients obtained independently is of significantly much better quality than the one using only one expansion.

\subsection{Nozzle flow}\label{Nozzle flow}
\begin{figure}[h!]
\subfigure[Mach]{\includegraphics[width=0.45\textwidth]{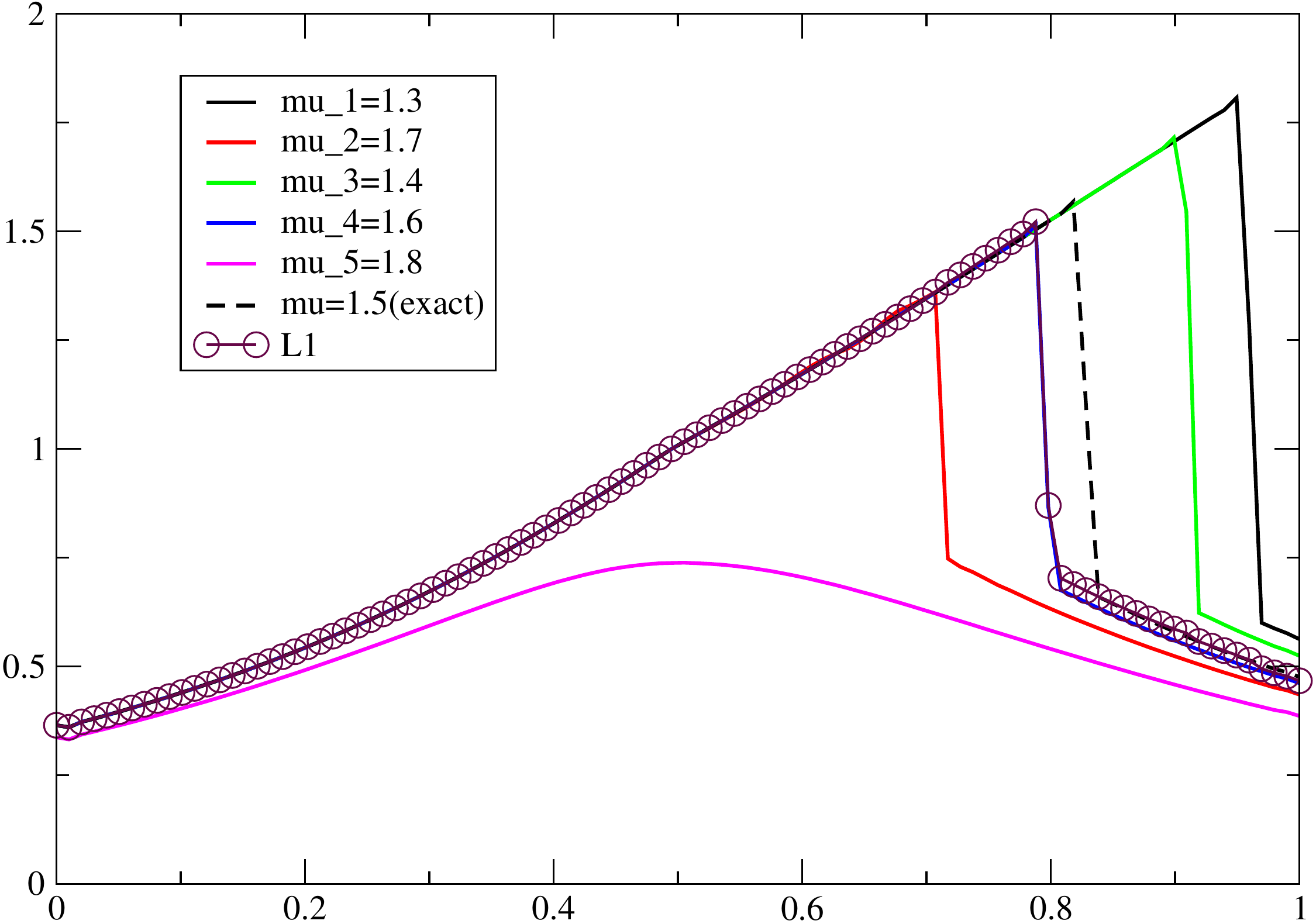}}
\subfigure[Density]{\includegraphics[width=0.45\textwidth]{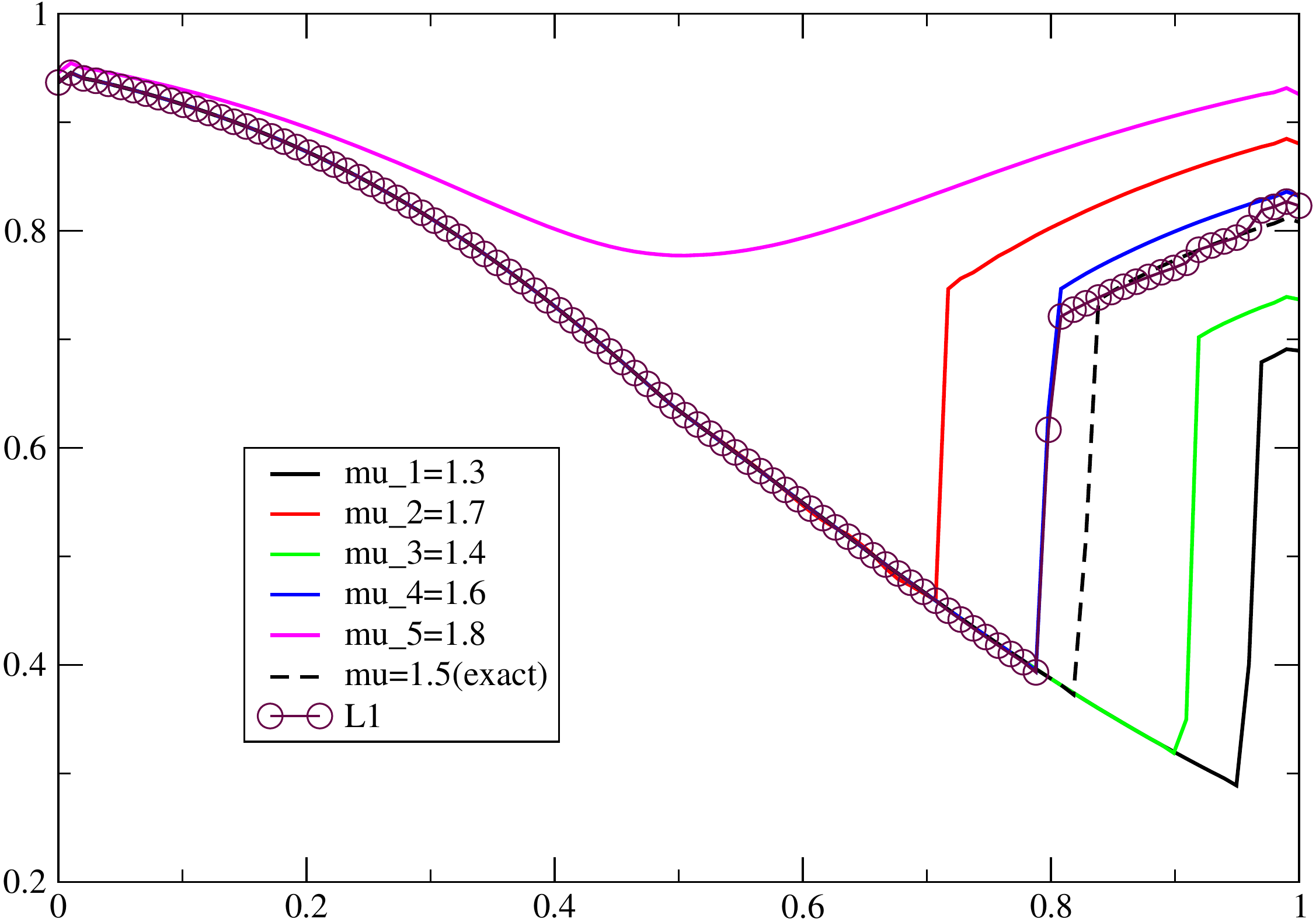}}
\caption{\label{Nozzle} Steady Nozzle flow: predicted solutions at target parameter $\mu=1.5$}
\end{figure}

To illustrate the ability of the reduced model, we consider the nozzle flow numerical experiment. The PDF is $$\frac{\partial F}{\partial x}=S(U)$$ where 
$$U=(\rho,\rho u, E)^{T},F(U)=(A\rho u,A(\rho u^{2}+p),Au(E+p))^{T}, S(U)=(0,p\frac{\partial A}{\partial x},0)^{T}$$ and $A$ is the area of the nozzle flow. Depending on the boundary conditions, we can have a fully smooth flow or a flow with steady discontinuity. We illustrate the method on a case where a discontinuity exists (see Figure \ref{Nozzle}). All the other variables behave in the same manner. The experiment has been conducted for the density case with the choice of the target parameter $\mu=1.5$. The nozzle is a Laval nozzle with
area given, for $x\in [0,1]$, by
$$A(x)=\left \{\begin{array}{ll}
1+6\big (x-\frac{1}{2}\big )^2 & \text{ if } 0\leq x\le \frac{1}{2}\\
1+0.15 \big (x-\frac{1}{2}\big )^2+6\big (x-\frac{1}{2}\big )^3& \text{ if } \frac{1}{2} \leq x\le 1
\end{array}
\right .
$$

The $L^1$ solution is the element of the dictionary which is the closest from the solution, as expected.

\section{Conclusions}
We have presented a general framework to approximate the solution of steady and unsteady hyperbolic problems. The solution can be smooth or discontinuous, and in the unsteady case, moving wave may exist. 
Starting from any  standard scheme (explicit or implicit), the reduced order solution is obtained at each time step (or each iteration in the implicit case) by from a minimization problem in the $L^1$ norm. We give a sufficient condition to be able to solve the problem, and discuss the practical aspects of the method. It is illustrated by several examples dealing with linear and non linear problems, scalar and systems in one space dimension. A rough error estimate is given, relating the successive projections and the initial solution.

This work can be extended in several directions. First, its efficiency can be improved by using hyper-reduction. The present contribution should be more seen as a proof of concept than a working algorithm, because its cost is still high. Hyper-reduction is probably the way to reduce drastically its cost. The second direction is to tackle multidimensional problems, and this can only be possible if the cost issue is solved. We are currently working on these two aspects.

\section*{Acknowledgements}
The first author was funded in part by the
MECASIF project  (2013-2017) funded by the French "Fonds Unique Interministériel" and SNF grant \# 200021\_153604 of the Swiss National Foundation.
The second  author would like to acknowledge partial support by the Army Research Laboratory through the Army High Performance Computing Research Center under Cooperative Agreement W911NF- 07-2-0027, and partial support by the Office of Naval Research under
grant no. N00014-11-1-0707. This document does not necessarily reflect the position of these institutions, and no official endorsement
should be inferred.
The third author has been supported in part by the SNF grant \# 200021\_153604 of the Swiss National Foundation

RA would also like for several very insightfull conversations with Y. Maday, LJLL, Universit\'e Pierre et Marie Curie.

\bibliographystyle{siam}

\begin{thebibliography}{10}

\bibitem{amsallem10}
{\sc D~Amsallem, J~Cortial, and C~Farhat}, {\em {Toward real-time
  computational-fluid-dynamics-based aeroelastic computations using a database
  of reduced-order information}}, AIAA Journal, 48 (2010), pp.~2029--2037.

\bibitem{amsallem13:mpc}
{\sc D~Amsallem, S~Deolalikar, F~Gurrola, and C~Farhat}, {\em {Model Predictive
  Control under Coupled Fluid-Structure Constraints Using a Database of
  Reduced-Order Models on a Tablet}}, AIAA Paper 2013-2588, 21st AIAA
  Computational Fluid Dynamics Conference, San Diego, CA, June 26-29, 2013,
  (2013), pp.~1--12.

\bibitem{amsallem08}
{\sc D~Amsallem and C~Farhat}, {\em {Interpolation method for adapting
  reduced-order models and application to aeroelasticity}}, AIAA Journal, 46
  (2008), pp.~1803--1813.

\bibitem{amsallem15:smo}
{\sc D~Amsallem, M~J Zahr, Y~Choi, and C~Farhat}, {\em {Design Optimization
  Using Hyper-Reduced-Order Models}}, Structural and Multidisciplinary
  Optimization, 51 (2015), pp.~919--940.

\bibitem{amsallem12:localROB}
{\sc D~Amsallem, M~J Zahr, and C~Farhat}, {\em {Nonlinear model order reduction
  based on local reduced-order bases}}, International Journal for Numerical
  Methods in Engineering, 92 (2012), pp.~891--916.

\bibitem{balajewicz15}
{\sc M~Balajewicz, D~Amsallem, and C~Farhat}, {\em {Projection-based model
  reduction for contact problems}}, arXiv.org,  (2015).

\bibitem{barone09}
{\sc M~F Barone, I~Kalashnikova, DJ~Segalman, and HK~Thornquist}, {\em {Stable
  Galerkin reduced order models for linearized compressible flow}}, Journal of
  Computational Physics, 228 (2009), pp.~1932--1946.

\bibitem{binev11}
{\sc P~Binev, A~Cohen, W~Dahmen, R~DeVore, G~Petrova, and P~Wojtaszczyk}, {\em
  {Convergence Rates for Greedy Algorithms in Reduced Basis Methods}}, SIAM
  Journal on Mathematical Analysis, 43 (2011), pp.~1457--1472.

\bibitem{boyd04}
{\sc S~Boyd and L~Vandenberghe}, {\em {Convex optimization}}, Cambridge
  university press, 2004.

\bibitem{brunton13}
{\sc S~L Brunton, J~H Tu, I~Bright, and J~N Kutz}, {\em {Compressive sensing
  and low-rank libraries for classification of bifurcation regimes in nonlinear
  dynamical systems}}, arXiv.org,  (2013).

\bibitem{buithanh08}
{\sc T~Bui-Thanh, K~Willcox, and O~Ghattas}, {\em {Parametric reduced-order
  models for probabilistic analysis of unsteady aerodynamic applications}},
  AIAA Journal, 46 (2008), pp.~2520--2529.

\bibitem{candes06}
{\sc E~Candes and J~Romberg}, {\em {Robust Signal Recovery from Incomplete
  Observations}}, in Image Processing, 2006 IEEE International Conference on,
  IEEE, 2006, pp.~1281--1284.

\bibitem{carlberg11}
{\sc K~Carlberg, C~Bou-Mosleh, and C~Farhat}, {\em {Efficient non-linear model
  reduction via a least-squares Petrov--Galerkin projection and compressive
  tensor approximations}}, International Journal for Numerical Methods in
  Engineering, 86 (2011), pp.~155--181.

\bibitem{carlberg13}
{\sc K~Carlberg, C~Farhat, J~Cortial, and D~Amsallem}, {\em {The GNAT method
  for nonlinear model reduction: effective implementation and application to
  computational fluid dynamics and turbulent flows}}, Journal of Computational
  Physics, 242 (2013), pp.~623--647.

\bibitem{chaturantabut10}
{\sc S~Chaturantabut and DC~Sorensen}, {\em {Nonlinear model reduction via
  discrete empirical interpolation}}, SIAM Journal on Scientific Computing, 32
  (2010), pp.~2737--2764.

\bibitem{dahmen14}
{\sc W~Dahmen, C~Plesken, and G~Welper}, {\em {Double greedy algorithms:
  Reduced basis methods for transport dominated problems}}, ESAIM: Mathematical
  Modelling and Numerical Analysis, 48 (2014), pp.~623--663.

\bibitem{daubechies08}
{\sc I~Daubechies, R~DeVore, Massimo Fornasier, and C~Sinan Gunturk}, {\em
  {Iteratively re-weighted least squares minimization for sparse recovery}},
  arXiv.org,  (2008).

\bibitem{dihlmann11}
{\sc M~Dihlmann, M~Drohmann, and B~Haasdonk}, {\em {Model reduction of
  parametrized evolution problems using the reduced basis method with adaptive
  time-partitioning}}, Proc. of ADMOS, 2011 (2011).

\bibitem{donoho06}
{\sc D~L Donoho}, {\em {Compressed sensing}}, IEEE Transactions on Information
  Theory, 52 (2006), pp.~1289--1306.

\bibitem{God}
{\sc E~Godlewski and PA~Raviart}, {\em {Hyperbolic systems of conservation
  laws}}, Ellipses, Feb. 1991.

\bibitem{grepl05}
{\sc M~A Grepl and A~T Patera}, {\em {A posteriori error bounds for
  reduced-basis approximations of parametrized parabolic partial differential
  equations}}, ESAIM: Mathematical Modelling and Numerical Analysis, 39 (2005),
  pp.~157--181.

\bibitem{guermond08}
{\sc J~L Guermond, F~Marpeau, and B~Popov}, {\em {A fast algorithm for solving
  first-order PDEs by L1-minimization}}, Communications in Mathematical
  Sciences, 6 (2008), pp.~199--216.

\bibitem{guermond09}
{\sc J~L Guermond and B~Popov}, {\em {$L^1$-Approximation of Stationary
  Hamilton--Jacobi Equations}}, SIAM Journal on Numerical Analysis, 47 (2009),
  pp.~339--362.

\bibitem{hovland08}
{\sc S.~Hovland, J.T. Gravdahl, and K~Willcox}, {\em {Explicit model predictive
  control for large-scale systems via model reduction}}, Journal of Guidance,
  Control, and Dynamics, 31 (2008), pp.~1--23.

\bibitem{huber11}
{\sc Peter~J Huber and Elvezio~M Ronchetti}, {\em {Robust Statistics}}, John
  Wiley {\&} Sons, Sept. 2011.

\bibitem{ito98}
{\sc K~Ito and SS~Ravindran}, {\em {A Reduced-Order Method for Simulation and
  Control of Fluid Flows}}, Journal of Computational Physics, 143 (1998),
  pp.~403--425.

\bibitem{kaulmann13}
{\sc S~Kaulmann and B~Haasdonk}, {\em {Online Greedy Reduced Basis Construction
  Using Dictionaries}}, University of Stuttgart,  (2013).

\bibitem{kunisch02}
{\sc K~Kunisch}, {\em {Galerkin proper orthogonal decomposition methods for a
  general equation in fluid dynamics}}, SIAM Journal on Numerical Analysis,
  (2003).

\bibitem{kunisch01}
{\sc K~Kunisch and S~Volkwein}, {\em {Galerkin proper orthogonal decomposition
  methods for parabolic problems}}, Numerische Mathematik, 90 (2001),
  pp.~117--148.

\bibitem{lee99}
{\sc D~D Lee and H~S Seung}, {\em {Learning the parts of objects by
  non-negative matrix factorization}}, Nature, 401 (1999), pp.~788--791.

\bibitem{legresley00}
{\sc P~A LeGresley and J~J Alonso}, {\em {Airfoil design optimization using
  reduced order models based on proper orthogonal decomposition}}, AIAA Paper
  2000-2545 Fluids 2000 Conference and Exhibit, Denver, CO,  (2000), pp.~1--14.

\bibitem{ly01}
{\sc HV~Ly and HT~Tran}, {\em {Modeling and control of physical processes using
  proper orthogonal decomposition}}, Mathematical and computer modelling, 33
  (2001), pp.~223--236.

\bibitem{maday12}
{\sc Y~Maday and B~Stamm}, {\em {Locally adaptive greedy approximations for
  anisotropic parameter reduced basis spaces}}, arXiv.org,  (2012).

\bibitem{nocedalbook}
{\sc J~Nocedal and S~J Wright}, {\em {Numerical optimization}}, Springer, Dec.
  2006.

\bibitem{pdt15}
{\sc A~Paul-Dubois-Taine and D~Amsallem}, {\em {An adaptive and efficient
  greedy procedure for the optimal training of parametric reduced-order
  models}}, International Journal for Numerical Methods in Engineering, 102
  (2015), pp.~1262--1292.

\bibitem{rozza08}
{\sc G~Rozza, DBP Huynh, and A~T Patera}, {\em {Reduced Basis Approximation and
  a Posteriori Error Estimation for Affinely Parametrized Elliptic Coercive
  Partial Differential Equations}}, Archives of Computational Methods in
  Engineering, 15 (2008), pp.~229--275.

\bibitem{ryckelynck05}
{\sc D~Ryckelynck}, {\em {A priori hyperreduction method: an adaptive
  approach}}, Journal of Computational Physics, 202 (2005), pp.~346--366.

\bibitem{sirovich87}
{\sc L~Sirovich}, {\em {Turbulence and the dynamics of coherent structures.
  Part I: coherent structures}}, Quarterly of applied mathematics, 45 (1987),
  pp.~561--571.

\bibitem{toro}
{\sc E.F. Toro}, {\em Riemann solvers and numerical methods for fluid
  dynamics}, Springer, Berlin, Heidelberg, 1997.

\bibitem{veroy03}
{\sc K~Veroy, C~Prud'homme, D~V Rovas, and A~T Patera}, {\em {A posteriori
  error bounds for reduced-basis approximation of parametrized noncoercive and
  nonlinear elliptic partial differential equations}}, AIAA Paper 2003-3847,
  16th Computational Fluid Dynamics Conference, Orlando Florida, 24-25 June
  2003, 3847 (2003).

\bibitem{willcox01}
{\sc K~Willcox and J~Peraire}, {\em {Balanced model reduction via the proper
  orthogonal decomposition}}, AIAA 2001-2611,  (2001), pp.~1--9.

\end{thebibliography}

\end{document}